\documentclass[11pt,a4paper]{article}%
\usepackage{amsmath}
\usepackage{graphicx}
\usepackage{epsfig}%
\usepackage{amsfonts}%
\usepackage{amssymb}
\usepackage{algorithm}
\usepackage{algorithmic}
\usepackage{color}

\makeatletter\@addtoreset{equation}{section}\makeatother
\makeatletter\@addtoreset{figure}{section}\makeatother
\makeatletter\@addtoreset{table}{section}\makeatother
\begin{document}

\title{Tensor Numerical Methods in Quantum Chemistry: \\
from Hartree-Fock Energy  to Excited States}

\author{V. KHOROMSKAIA\thanks{Max-Planck-Institute for
        Mathematics in the Sciences, Leipzig; Max-Planck-Institute for Dynamics of Complex Systems, Magdeburg,
        Germany ({\tt vekh@mis.mpg.de}).}
        \and B. N. KHOROMSKIJ\thanks{Max-Planck-Institute for
        Mathematics in the Sciences, Inselstr.~22-26, D-04103 Leipzig,
        Germany ({\tt bokh@mis.mpg.de}).}}

 \date{}

 \maketitle

\begin{abstract}
We  resume the recent successes of the grid-based tensor numerical methods and 
discuss their prospects in   real-space electronic structure calculations.
These methods, based on the low-rank representation of the multidimensional functions and 
integral  operators, first appeared as an accurate tensor calculus for the 3D Hartree potential
using 1D complexity operations, and have evolved  to
entirely grid-based tensor-structured 3D Hartree-Fock eigenvalue solver.
It  benefits from  tensor calculation of the core Hamiltonian
and two-electron integrals (TEI) in $O(n\log n)$ complexity using the
rank-structured  approximation of basis functions, electron densities and  convolution 
integral  operators  all represented  on 3D $n\times n\times n $ Cartesian grids.
The algorithm for calculating TEI tensor in a form of the Cholesky decomposition
is based on multiple factorizations using  algebraic 1D ``density fitting``  scheme , 
which yield an almost irreducible number of product basis functions  involved in the
$3$D convolution integrals,  depending on a threshold $\varepsilon >0$. 
The basis functions are not restricted to separable Gaussians, since the analytical integration is 
 substituted by high-precision tensor-structured numerical quadratures. 
The tensor approaches to post-Hartree-Fock calculations for the MP2  energy  
correction and for the Bethe-Salpeter excited states, based on  using   
low-rank factorizations and the reduced basis method,  were recently introduced.
 Another direction is related to the recent attempts to develop a tensor-based  
Hartree-Fock numerical scheme  for finite lattice-structured systems, where one of
the numerical challenges is the summation of electrostatic potentials of a large number of nuclei.
The  3D grid-based  tensor method for calculation of a potential sum 
on a $L\times L\times L$ lattice  manifests the linear in $L$ computational work,  
$O(L)$,  instead of the usual $O(L^3 \log L)$ scaling by the Ewald-type approaches. 
The accuracy of the order of atomic radii, $ h\sim 10^{-4} \AA{}$,  for the grid 
representation of electrostatic potentials  is achieved due to low cost of using  $1$D 
operations on  large 3D grids.
\end{abstract}

\maketitle



 \section{Introduction}

The  problems of numerical modeling the many-particle interactions in large molecular
systems, lattice structured metallic clusters and crystals, proteins and nanomaterials
are the most challenging tasks in modern computational physics and chemistry.
The traditional approaches for these multidimensional problems are restricted to
the concepts which have well recognized limitations for computations with higher accuracy and for 
larger non-periodic molecular systems, as well as for efficient calculation of excited states.
Recent advances in numerical analysis for multidimensional problems and significant achievements 
in high performance computing suggest new creative approaches to these problems.

The Hartree-Fock (HF) equation governed by the 3D integral-differential operator 
\cite{SzOst:1996,HJO:2000} is one of the basic models for \emph{ab initio} 
calculation of the ground state energy of molecular systems.
It is a strongly nonlinear eigenvalue problem (EVP) in a sense, that one should  
solve the equation in a self-consistent way when the integral part of the governing 
operator  depends on the solution itself.
Multiple strong singularities in the electron density 
of a molecule due to nuclear cusps impose  strong requirements on the accuracy of calculations.

Commonly used numerical methods for solution of the Hartree-Fock equation are based on the analytical 
computation of the arising two-electron integrals (convolution type integrals in $\mathbb{R}^3$)
in the problem adapted naturally separable Gaussian-type bases \cite{Alml:95},
by using {\it erf}-function  expansions. 
This rigorous approach resulted in a  number of efficient implementations 
which are widely used in computational quantum chemistry. The success of the analytical
integration methods stems from the big amount of precomputed information based on the physical
 insight  including the construction of problem adapted atomic orbitals basis sets and 
 elaborate nonlinear optimization for calculation of density fitting basis.
 The known  limitations of this approach appear due to a strong dependence of 
the numerical efficiency on the size and quality of the chosen Gaussian basis sets, 
 that might be crucial  for larger molecular clusters and heavier atoms.

The intention to replace or assist the analytical  calculations 
for the Hartree-Fock problem by a data-sparse 
grid-based  numerical schemes  has a long history. First success was the grid-based numerical method for 
the diatomic molecules in \cite{SuPyLa:1985}, though this approach was not feasible to compact (3D) molecules. 
The wavelet multiresolution schemes \cite{DamPrSch:95} capable 
for the accurate representation of nuclear cusps,
 have been applied to electronic structure simulations since the seminal papers 
\cite{HaFaYaBeyl:04,SeMaYaHa:08}, and recently this approach was further advanced 
due to achievements in the high performance supercomputing 
\cite{SeMaYaHa:08,BiHaVa:12,FreFoFlaRu:13,DuJenJuWiFlaFre:14}.
However,  due to extensive computational resources, 
the entirely  wavelet-based or sparse-grid approaches \cite{Griebel:07,Yseren:10} 
are limited so far only to rather small atomic systems with few electrons \cite{Bach:12}.
Tensor hyper-contraction decompositions in density fitting schemes have been analyzed 
in \cite{HoPaMa:12,PaHoMaSche:13}.

The newly developed tensor-structured  numerical methods, both the name and the concept,
appeared during the work on the grid-based tensor approach to the solution of the 3D
Hartree-Fock problem \cite{KhKh:06,KhKh3:08,KhKhFl_Hart:09}. 
The central point is the representation of $d$-variate functions and 
integral/differential operators on large $n^{\otimes d }$ grids and their approximation  
in the  low-rank tensor  formats, which
allows  numerical calculations in $O(dn)$ complexity instead of $O(n^d)$ by 
conventional methods.

In this paper we  summarize the main benefits of the  
tensor numerical methods in electronic structure calculations, 
and discuss further prospects of this approach  
in several directions, such as
\begin{itemize}
 \item algebraic directional density fitting and tensor factorization of the two-electron integrals 
and their use in the Hartree-Fock calculations;

\item tensor decompositions in the MP2 energy correction, and in calculation of excited states
based on the Bethe-Salpeter equation;

\item fast tensor summation of electrostatic potentials on large 3D lattice for
efficient calculation of the Fock matrix in the case of lattice-structured systems;

\item fast calculation of the interaction energy of the long-range potentials 
on large 3D lattices.
\end{itemize} 
Notice that basic rank-structured tensor formats such as the canonical (PARAFAC) 
and Tucker tensor decompositions have been since
long used in the computer science for the quantitative analysis of correlations in the experimental
multidimensional data arrays in data processing and chemometrics, 
see \cite{KoldaB:08} and references therein.
In 2006 the exceptional properties of the Tucker decomposition for the discretized
multidimensional functions have been revealed in \cite{Khor:06,KhKh:06},
where it was proven that for a class of function-related tensors the approximation error of the
Tucker decomposition {\it decays exponentially} with respect to the Tucker rank.
This gives the opportunity to represent the discretized multidimensional functions
and integral (convolution) operators
in an algebraically separable form and thus reduce the numerical treatment of the
multidimensional  transforms  to 1D operations. It was shown that the number of vectors in such a
representation  depends  only logarithmically on the size of the $d$-dimensional grid,
$O(n^d)$, used for discretization of multivariate functions.

The above results have led to the idea to calculate the Hartree potential and the Coulomb and
exchange operators by numerical quadratures \cite{KhKh3:08} using  
Gaussian-type basis functions discretized on 3D Cartesian grids and the 
fast tensor-product convolution \cite{Khor1:08}. 
In this way, the efficient low-rank canonical tensor representation to the Newton 
kernel was an essential contribution \cite{BeHaKh:08}. 
To reduce the initial rank of the electron density,
that is quadratically proportional to  the number of GTO basis functions, 
the canonical-to-Tucker tensor transform was invented \cite{KhKh3:08,KhKh:06}, 
which made computations for large tensor grids and extended molecules tractable (even in Matlab).

The initial version of tensor-structured algorithms for solving Hartree-Fock equation 
employed the 3D grid-based calculation of the Coulomb and exchange integral 
operators ''on-the-fly'', thus avoiding precomputation and storage of the TEI 
tensor\cite{VKH:exch2010,KhKhFl_Hart:09}.
In particular, it was justified that tensor calculus allows to reduce 
the 3D convolution integrals to combinations of 1D convolutions, and 1D Hadamard 
and scalar products.
Besides, these results promoted spreading of the tensor-structured methods 
in the community of numerical analysis 
\cite{GraKresTo:13,KhorSurv:14,HaSchneid_Surv:15},
and further development of the tree-tensor formats like tensor-train \cite{Osel_TT:11} 
and hierarchical Tucker representations \cite{HaKu:08}.

Further development of tensor methods in electronic structure calculations 
was due to the fast algorithm for the grid-based computation of the 
TEI tensor \cite{VeKhBoKhSchn:12} 
in $O(N_b^3)$ storage in the number of basis functions $N_b$. 
The fourth order TEI tensor is calculated in a form of the Cholesky factorization by using
the algebraic 1D ''density fitting`` scheme, which applies to the 
product basis functions. 
Imposing the low-rank tensor representation of the product basis functions and 
the Newton convolving kernel all discretized on $n\times n\times n $ Cartesian grid,
the 3D integral transforms are calculated in $O(n\log n)$ complexity. 
Given the factorized TEI, the update of the Column and exchange parts in the Fock matrix 
reduces to the cheap algebraic operations.
Other steps are tensor calculation of the core Hamiltonian and the efficient 
MP2 energy correction scheme \cite{VeBoKhMP2:13},
which all together gave rise to the black-box Hartree-Fock solver \cite{VeKh_box:13}.

Due to the grid representation of basis functions  basis sets are now not restricted to 
Gaussian-type orbitals and allowed to be any well-separable function defined on a grid. 
The tensor-based Hartree-Fock solver  is competitive in computational time and accuracy 
with the standard packages based on analytical calculation of integrals.
High accuracy is attained owing to easy calculations on 
large 3D grids up to $n^3=10^{18}$,
so that the resolution with mesh size $h$ of the order of atomic radii, 
$10^{-4} \AA{}$, is possible.

Motivated by tensor decompositions in the MP2 scheme, 
the new approach for calculation of the excited states  in the framework of the 
Bethe-Salpeter equation was recently introduced \cite{BeKhKh_BSE:15}, 
that employs  the reduced basis method in combination with low-rank tensor approximations.

Further developments of the tensor methods in multi-particle simulations are focused on the 
large lattice structures in a box and nearly periodic systems, for which the  
tensor approach may reduce computational costs dramatically \cite{KhKh_arx:14}. 
One of the challenging problems
in the numerical treatment of the crystalline-type molecular clusters is summation of a large number
of electrostatic potentials distributed on a finite (non-periodic) lattice.
The novel grid-based method for summation of the long-range potentials in the canonical 
and Tucker formats \cite{KhKh_CPC:13,VeBoKh_Tuck:14}
works on $L\times L\times L $ lattices with the computational cost  $O(L)$ 
instead of $O(L^3 \log L)$ by the traditional Ewald-type methods \cite{Ewald:27}. 
The required precision is guaranteed  by employing large 3D Cartesian grids 
for representation of potentials. 
The method remains efficient for lattices with non-rectangular 
geometries and in the presence of multiple defects \cite{VeBoKh_Tuck:14,KhKh_arx:14}.

The rest of the paper is organized as follows. 
In Section \ref{sec:Tens_formats}, we overview the basic tensor
formats and show why the orthogonal Tucker tensor decomposition, originating from computer science
and data processing, became useful for the treatment of the multidimensional
functions and operators in numerical analysis. In particular, \S \ref{ssec:MPS-TT_formats} and
\S \ref{ssec:QQTform} discuss the matrix product states (tensor train) formats 
and the novel quantics tensor approximation method, respectively, while
\S \ref{ssec:Oper_1D} addresses the basic tensor-structured multilinear algebra operations.
Section \ref{sec:Tensor_HF} describes tensor calculus of the multidimensional
convolution transform on examples of the Hartree potential (\S \ref{ssec:Multi_Int}) and the 
TEI tensor (\S \ref{sssec:TEI_Int}), and recalls the main building blocks 
in the tensor-based black-box Hartree-Fock solver.
Section \ref{sec:MP2_BSE} shows benefits of the tensor approach in MP2 calculations
and in calculation of the lowest part  of the excitation energies in the framework of the 
Bethe-Salpeter equation. Section \ref{sec:Cryst} describes the benefits of tensor methods
in applications to lattice-type molecular systems. In particular,
we present fast tensor method for the summation of the long-range interaction potentials on 
a 3D lattice by using the assembled vectors of their canonical and Tucker tensor representations.
The important application of this approach to calculation of interaction energy of the
Coulombic potentials on a lattice with sub-linear cost in the lattice size is described in detail.
Appendix discusses some computational details on the Canonical-to-Tucker tensor transform,
the Galerkin discretization scheme for the nonlinear Hartree-Fock equation, 
and the basics of the low-rank canonical tensor representation for the Newton kernel.

\section{Rank-structured tensor representations of discretized functions and operators}
\label{sec:Tens_formats}

In this section we discuss shortly why the rank-structured tensors, which were traditionally
used in experimental data processing, appears to be useful for the separable representation of 
multivariate functions and operators represented on tensor product grids. 
The canonical and Tucker tensor decompositions have been since 
long used in the computer science community for the quantitative analysis of correlations in the experimental 
multidimensional data arrays in chemometrics and data processing, with rather weak requirements
on the accuracy, see \cite{SmBroGe-Book:04,KoldaB:08}  and references therein.
In the recent decade the class of matrix product states  type tensor formats became popular in the
simulations of quantum spin systems and quantum molecular dynamics. 
We refer to \cite{KoldaB:08,Scholl:11,KhorSurv:14} for recent literature surveys
on commonly used tensor formats.

\subsection{Canonical and Tucker tensor formats}\label{ssec:CanTuck_formats}

We consider a tensor of order $d$, as a real  multidimensional array 
${\bf A}=[a_{i_1,...,i_d}] \in \mathbb{R}^{n_1 \times \ldots \times n_d}$
numbered by a $d$-tuple index 
set\footnote{The alternative notation ${\bf A}= [{\bf A}(i_1,...,i_d)]$ can be utilized.}, 
with multi-index notation ${\bf i}=(i_1,...,i_d) $, $i_\ell \in \{1,...,n_\ell\}, \ell=1,...,d$.
It is an element of the linear vector space $\mathbb{R}^{n_1 \times \ldots \times n_d}$ 
equipped with the Euclidean scalar product,
\[
 \left\langle {\bf A} ,{\bf B} \right\rangle = 
 \sum_{i_1=1}^{n_1} \cdots \sum_{i_d=1}^{n_d} a_{i_1,...,i_d} b_{i_1,...,i_d}.
\]
Euclidean vectors and matrices are the special case of $d$th order tensors.
For a general tensor, the  required storage 
scales exponentially in the dimension, $n_1 n_2 \cdots n_d$, 
(the so-called ''curse of dimensionality``). 
To get rid of exponential scaling in the dimension, one can apply  the 
rank-structured separable representations of multidimensional tensors.
The simplest separable element is given by rank-$1$ canonical tensor, 
\[ {\bf U} = {\bf u}^{(1)}\otimes ... 
 \otimes {\bf u}^{(d)}\in \mathbb{R}^{n_1 \times \ldots \times n_d},
\]
with entries
\[
 u_{i_1,\ldots, i_d}=  u^{(1)}_{i_1}\cdot \cdot\cdot u^{(d)}_{i_d},
\]
requiring only $n_1+ ... +n_d$ numbers to store it. 

A tensor in the $R$-term canonical format is defined by
 \begin{equation}\label{eqn:CP_form}
   {\bf U} = {\sum}_{k =1}^{R} c_k
   {\bf u}_k^{(1)}  \otimes \ldots \otimes {\bf u}_k^{(d)},
 \quad  c_k \in \mathbb{R},
\end{equation}
where ${\bf u}_k^{(\ell)}\in \mathbb{R}^{n_\ell}$ are normalized vectors, 
and $R$ is called the canonical rank of a tensor. It is convenient to introduce
 the so-called side matrices 
$U^{(\ell)}=[{\bf u}_1^{(\ell)}\ldots {\bf u}_R^{(\ell)}] \in \mathbb{R}^{n_\ell \times R}$, 
$\ell=1, 2,3$, obtained by concatenation of the canonical vectors 
${\bf u}_k^{(\ell)}$, $k=1,\ldots R$.
Now the storage cost is bounded by  $d R n$. For $d\geq 3$
computation of the canonical rank of a tensor ${\bf U}$, i.e.
the minimal number $R$ in representation (\ref{eqn:CP_form}) and the respective
decomposition, is an $N$-$P$ hard problem.
In the case $d=2$ the representation (\ref{eqn:CP_form}) is merely a rank-$R$ matrix.

\begin{figure}[tbh]
  \begin{center}
\includegraphics[width=6.0cm]{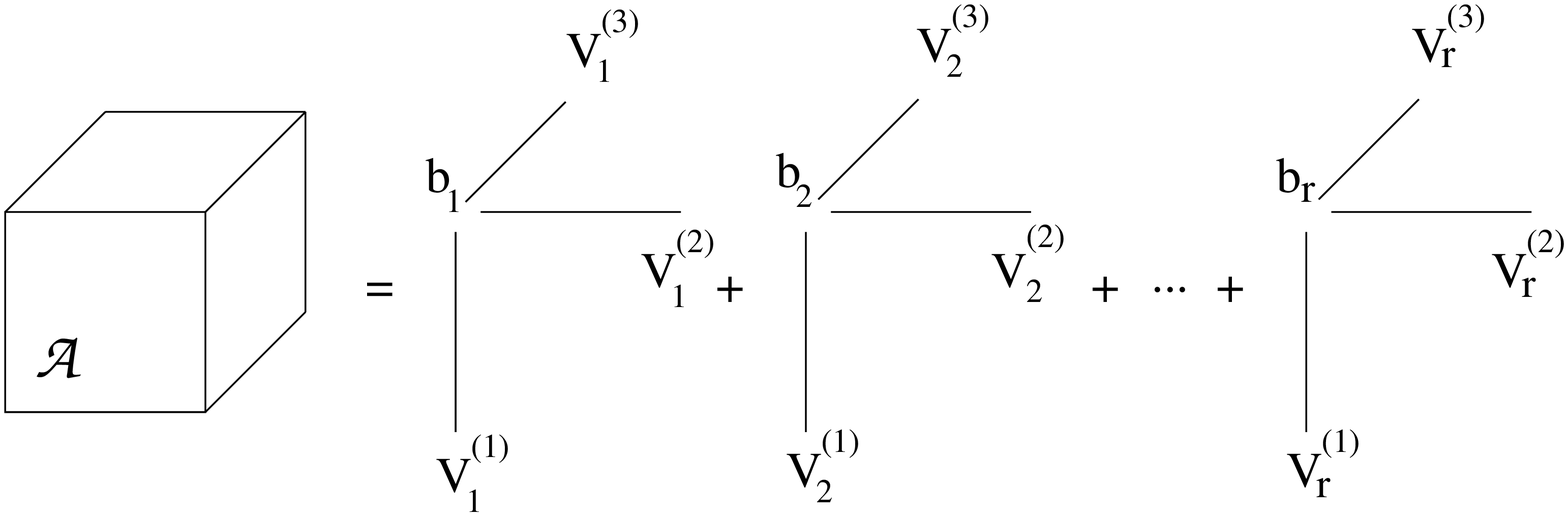} \quad 
\includegraphics[width=4.0cm]{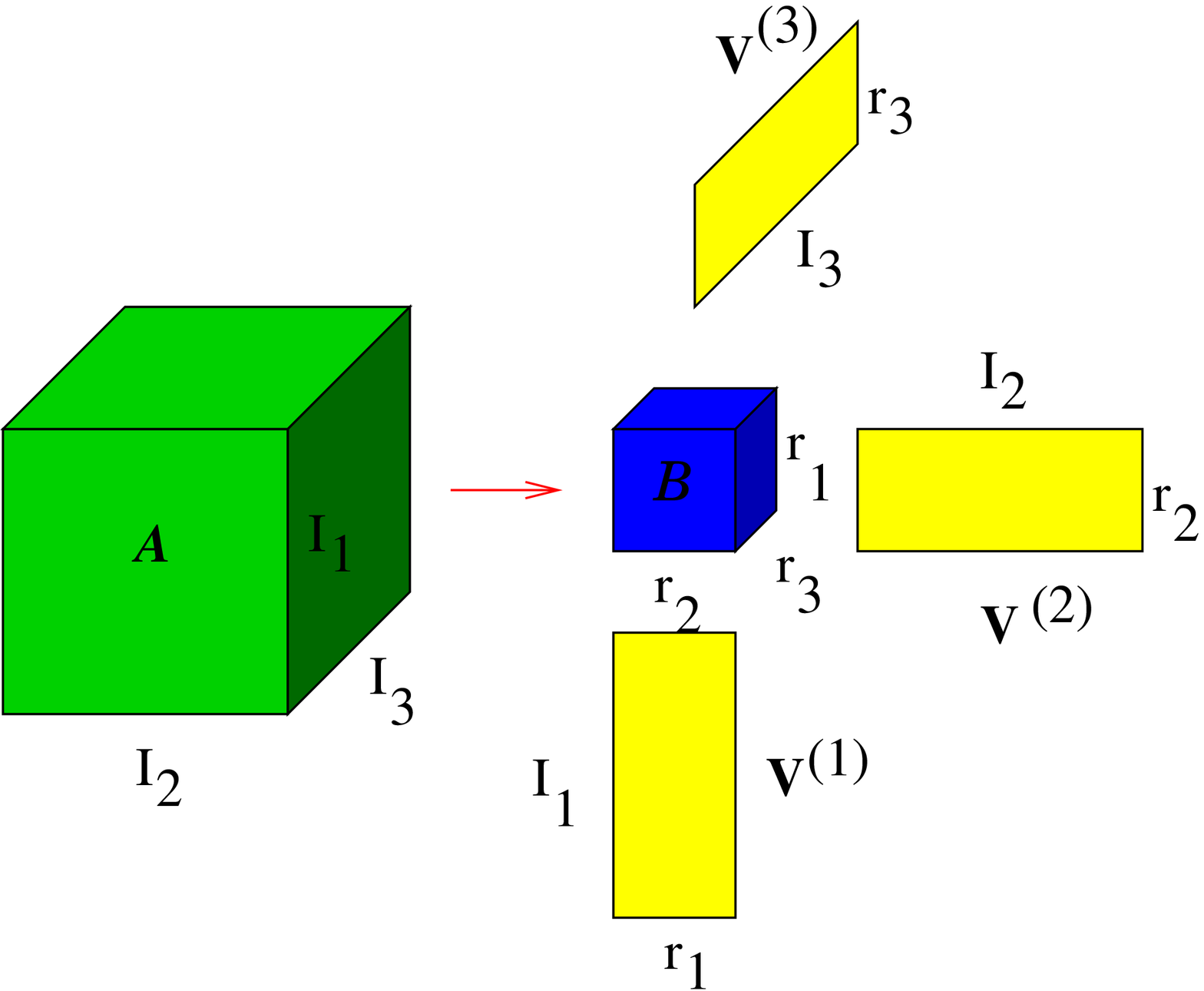}
\end{center}
\caption{Visualizing the canonical and  Tucker tensors for $d=3$.}
\label{fig:CanTucker_format}
\end{figure}

We say that a tensor ${\bf V}$ is represented in the rank-$\bf r$ orthogonal Tucker format 
with the rank parameter ${\bf r}=(r_1,...,r_d)$, if 
\begin{equation}\label{eqn:Tucker_form}
  {\bf V} =\sum\limits_{\nu_1 =1}^{r_1}\ldots
\sum\limits^{r_d}_{{\nu_d}=1} \beta_{\nu_1, \ldots ,\nu_d}
\,  {\bf v}^{(1)}_{\nu_1} \otimes \ldots \otimes {\bf v}^{(d)}_{\nu_d},\quad \ell=1,\ldots,d, 
\end{equation}
where $\{{\bf v}^{(\ell)}_{\nu_\ell}\in \mathbb{R}^{n_\ell}\}$ represents the set of orthonormal vectors, 
and $\boldsymbol{\beta}=[\beta_{\nu_1,...,\nu_d}] \in \mathbb{R}^{r_1\times \cdots \times r_d}$ is 
the Tucker core tensor. 
The storage cost for the Tucker tensor is bounded by $d r n +r^d$, $r=\max r_\ell$.
In the case $d=2$, the orthogonal Tucker decomposition is equivalent to the singular 
value decomposition (SVD) of a rectangular matrix. Figure \ref{fig:CanTucker_format}
visualizes the canonical and  Tucker tensors in the case $d=3$.


Canonical tensor decomposition leads to an ill-posed problem and, up to our best knowledge,  
there are no robust algebraic algorithms for the canonical approximation of an arbitrary tensor.
For some classes of analytic functions, the explicit low-rank approximation can be constructed
in analytic form by using $sinc$-quadrature approximation to the Laplace transform.
The robust methods for Tucker decomposition are based on the orthogonal projections
using the higher order singular value decomposition\cite{DMV-SIAM2:00} (HOSVD) 
that is the generalization of the singular value decomposition (SVD) of matrices.

The rank-structured decomposition (approximation) of multidimensional tensors provides
means for the separation of variables in the discretized representation of 
multivariate functions and operators,
and thus  the possibility to substitute multidimensional algebraic transforms 
by univariate operations. 
Notice that in the computer science community these possibilities were restricted 
to moderate-dimensional
tensors of small mode size (with rather large rank parameters) obtained 
from the experimental data sets.

\subsection{ Tucker decomposition for function related tensors. Canonical-to-Tucker approximation}
\label{ssec:Can2Tuck}

In 2006 it was  first verified numerically, see in  \cite{KhKh:06},  
that the rank of the (fixed accuracy) Tucker approximation to some function related tensors
depends only logarithmically on the size of the discretization grid.
For a given continuous function $f:\Omega\to \mathbb{R}$, 
$\Omega := \prod^d_{\ell=1} [a_\ell , b_\ell]\subset \mathbb{R}^3$, we 
introduced the function related $3$rd order tensor, obtained by Galerkin
discretization in a volume box using $n\times n\times  n$ 3D Cartesian grid.
In particular, the low-rank tensor approximations were calculated for functions 
of the Slater type $f(x)=\exp(-\alpha \|x\|)$,
Newton kernel $f(x)=\frac{1}{\| x \|}$, Yukawa potential $f(x)=\frac{e^{-\alpha \|x\| }}{\|x\|}$, 
and the Helmholtz potential $f(x)=\frac{\cos{\alpha \|x \| }}{\|x\|}$, $x \in \mathbb{R}^3 $. 

Numerical tests demonstrated that the error of the rank-${\bf r}$ (with ${\bf r}=(r,r,r)$) Tucker approximation  
applied to these third order tensors decays exponentially with respect to the Tucker rank $r$. 
Moreover, for fixed approximation error, the Tucker rank $r$  
scales as $O(\log n)$ for above mentioned functions \cite{Khor1:08}.
\begin{figure}[tbh]
\centering
 \includegraphics[width=4.5cm]{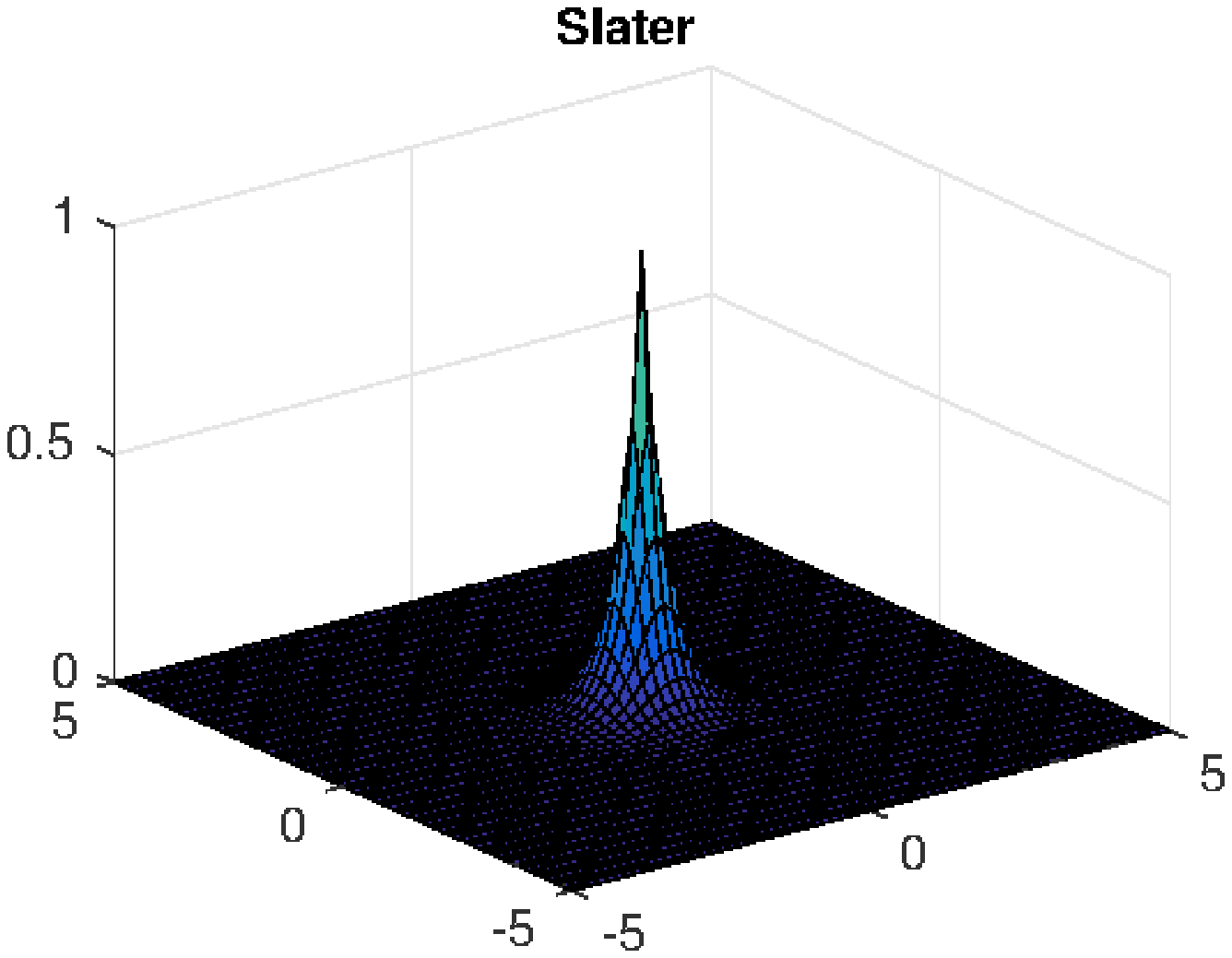}  \quad \quad
\includegraphics[width=4.6cm]{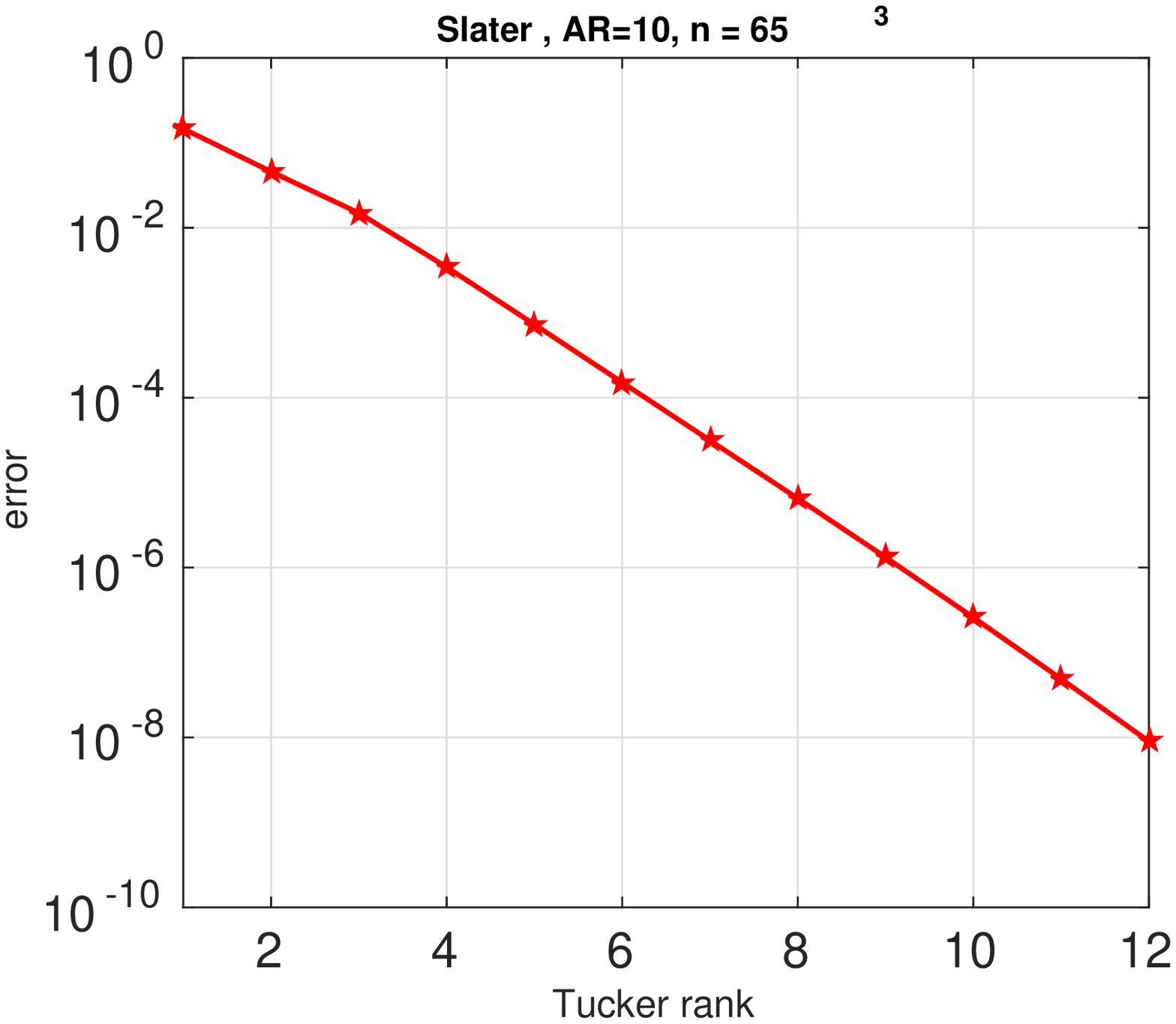}
\caption{The Tucker approximation error in the Frobenius norm vs. the Tucker rank $r$ for
a single  Slater function. }
\label{fig:Tucker_3D}
\end{figure}

This application revealed a number of properties of the Tucker format which were not 
known before. In particular, it was showed that for tensors resulting from the discretization 
of physically relevant multidimensional functions on the tensor grids one can find 
algebraically a reduced subspace in each mode of the original tensor, 
thus approximating the function by separated variables  using a 
tensor product of a relatively small number of vectors.

\begin{figure}[tbh]
\centering
\includegraphics[width=4.5cm]{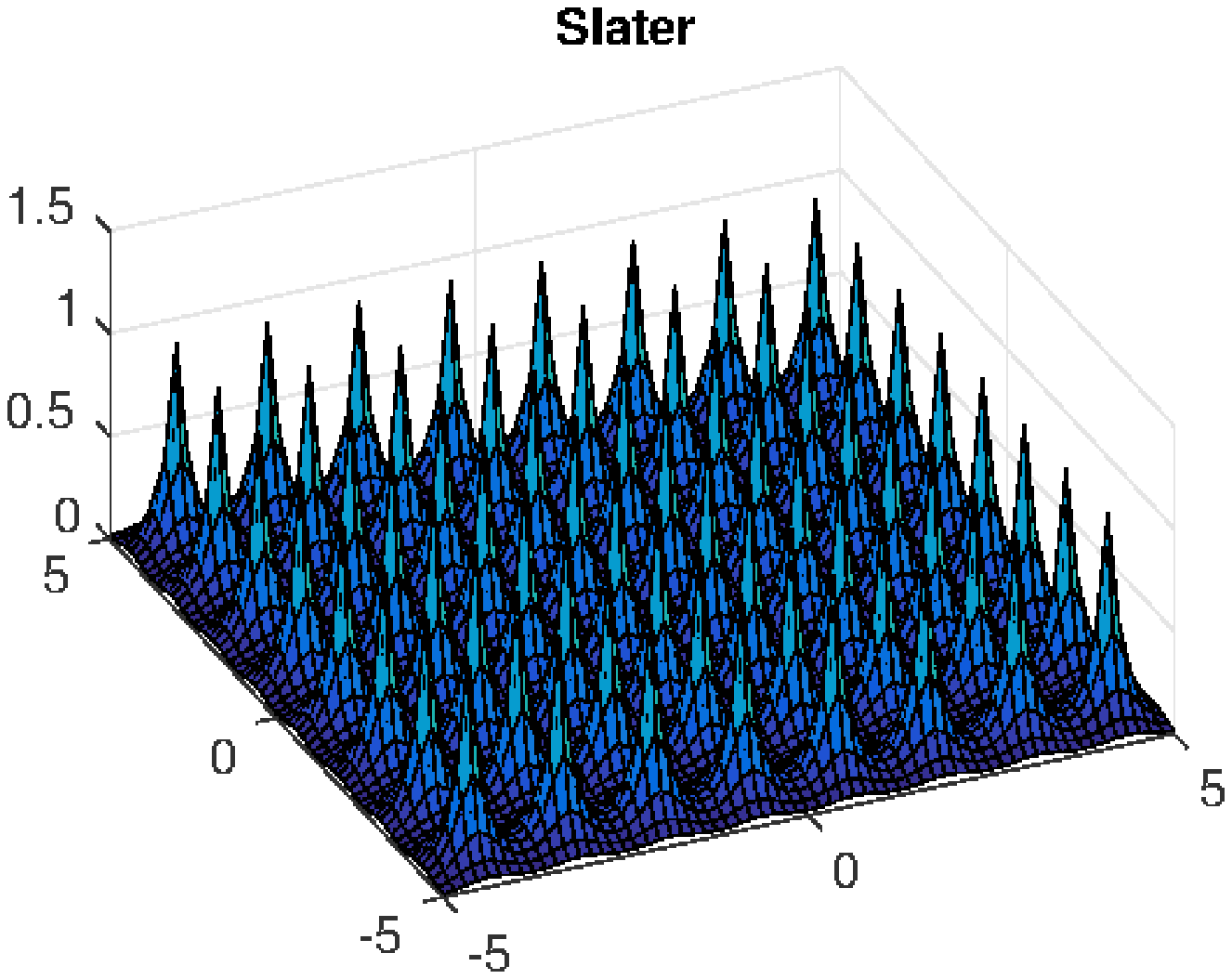} \quad \quad
\includegraphics[width=4.6cm]{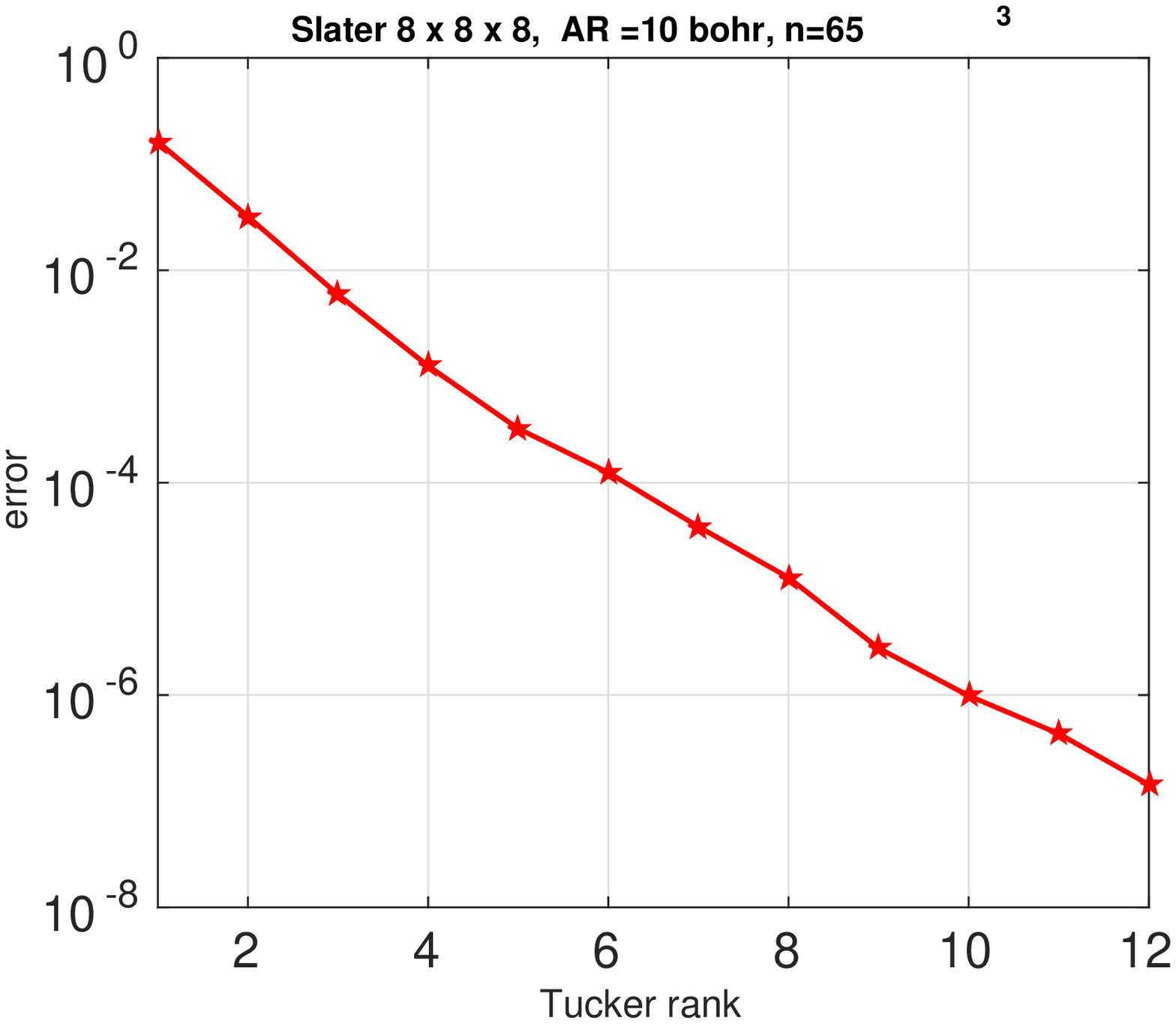} 
 \caption{The Tucker approximation error in the Frobenius norm vs. the Tucker rank $r$ for
a sum of Slater potentials over the cubic $8\times 8\times 8$ lattice (right). }
\label{fig:Tucker_3Dn}
\end{figure}

Figures \ref{fig:Tucker_3D} and \ref{fig:Tucker_3Dn} (right) show the exponential convergence 
of the Tucker tensor approximation in the Tucker rank $r$
of the Slater function in the relative Frobenius norm, 
$E_{FN}=\frac{\|A-A_r \|}{\|A\|}$.
Figures demonstrate that the exponential decay of approximation error is nearly the same 
for both a single potential  and for a sum of the same potentials distributed in nodes of a cubic lattice. 
This property  of the Tucker decomposition will be further gainfully applied to the fast 
lattice summation of interaction potentials.

Motivated by these observations, we invented the canonical-to-Tucker (C2T)
decomposition for function related tensors \cite{KhKh:06}, based on the reduced higher order
singular value decomposition (RHOSVD) (Theorem 2.5, \cite{KhKh3:08}).
The canonical-to-Tucker algorithm (see Appendix) combined with the Tucker-to-canonical
transform serves for reducing the ranks of canonical tensors with large $R$.
\begin{figure}[tbh]
\centering
   \includegraphics[width=10.5cm]{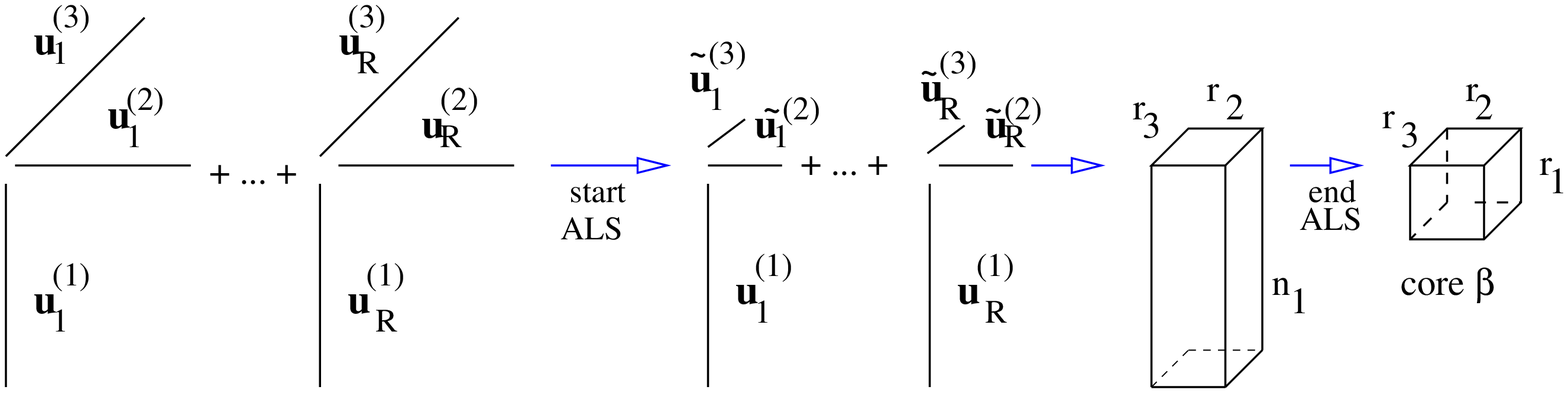}
\caption{Part of the canonical-to-Tucker decomposition algorithm for mode $\ell=1$.}
\label{fig:CanTucker_3D}
\end{figure}
Figure \ref{fig:CanTucker_3D} in its ALS part shows the algorithmic step for $\ell=1$, which is
repeated for every mode $\ell=1, 2, 3$ (see Appendix).
The computational work  for the multigrid tensor decomposition C2T algorithm  
introduced in \cite{KhKh3:08}
exhibits linear complexity scaling with respect to all input parameters,  $O(R n r)$.

The multigrid Tucker decomposition algorithm 
applied  to full format tensors \cite{KhKh3:08}
leads to the complexity  scaling $O(n^3)$,
whereas its standard
version (commonly used HOSVD algorithm in computer science) scales as $O(n^4)$,
becoming intractable even for moderate sizes of tensors.

We conclude by notice that 
the  optimized Tucker and canonical tensor approximations can be
computed by the alternating least square (ALS) iteration with initial guess 
obtained by HOSVD/RHOSVD approximations.

\subsection{Matrix Product States/Tensor Train format}\label{ssec:MPS-TT_formats}

The product-type representation of $d$th order tensors, which is
called in the physical literature as  the matrix product states (MPS) decomposition 
(or more generally, tensor network states models),
was introduced and successfully applied in DMRG quantum computations
\cite{White:93,Cirac_TC:04,Scholl:11}, and, independently, in quantum molecular dynamics 
as the multilayer (ML) MCTDH methods \cite{WangTho:03,LubBook:08}.
MPS type representations  reduce the storage complexity to $O(d r^2 N)$,
where $r$ is the maximal rank parameter. 

In the recent years the various versions of the MPS-type
tensor approximations were discussed and further investigated in mathematical literature 
including the hierarchical dimension splitting \cite{Khor:06}, 
the tensor train (TT) \cite{Osel_TT:11,DolgDiss:14}, the quantics-TT (QTT) \cite{KhQuant:09},
as well as the  hierarchical Tucker  representation \cite{HaKu:08}, which
belongs to the class of tensor network states model. 

The TT format is the particular form of MPS type factorization in the case of 
open boundary conditions. For a given rank parameter ${\bf r}=(r_0,...,r_d)$, 
the rank-${\bf r}$ TT format contains all elements
${\bf A}=[a_{i_1,...,i_d}] \in\mathbb{R}^{n_1\times ... \times n_d} $,
which can be represented as the contracted products of $3$-tensors, 
that in the index notation takes a form,
\[
a_{i_1,...,i_d}=\sum\limits_{\alpha_1 =1}^{r_1} \cdots  \sum\limits_{\alpha_d =1}^{r_d}
{\bf a}^{(1)}_{\alpha_1}(i_1)  {\bf a}^{(2)}_{\alpha_1,\alpha_2}(i_2)
\cdot\cdot\cdot {\bf a}^{(d)}_{\alpha_{d-1}}(i_d)
  \equiv {A}^{(1)}({i_1}) {A}^{(2)}({i_2}) \ldots {A}^{(d)}({i_d}),
\]
specified by the set of column vectors, ${\bf a}^{(\ell)}_{\alpha_\ell,\alpha_{\ell+1}}   \in \mathbb{R}^{n_\ell}$, 
($\ell=1,...,d$), or equivalently by the vector-valued $r_\ell\times r_{\ell+1}$ matrices,
${A}^{(\ell)}=[{\bf a}^{(\ell)}_{\alpha_\ell,\alpha_{\ell+1}}]$,
(i.e., $3$-tensors),  cf. (\ref{eqn:Tucker_form}).
The latter representation is written in the {\it matrix product form}, explaining the notion MPS, 
where ${A}^{(\ell)}(i_\ell) $ is $r_{\ell-1} \times r_\ell$ matrix.
\begin{figure}[tbh]
\begin{center}
\includegraphics[width=8.0cm]{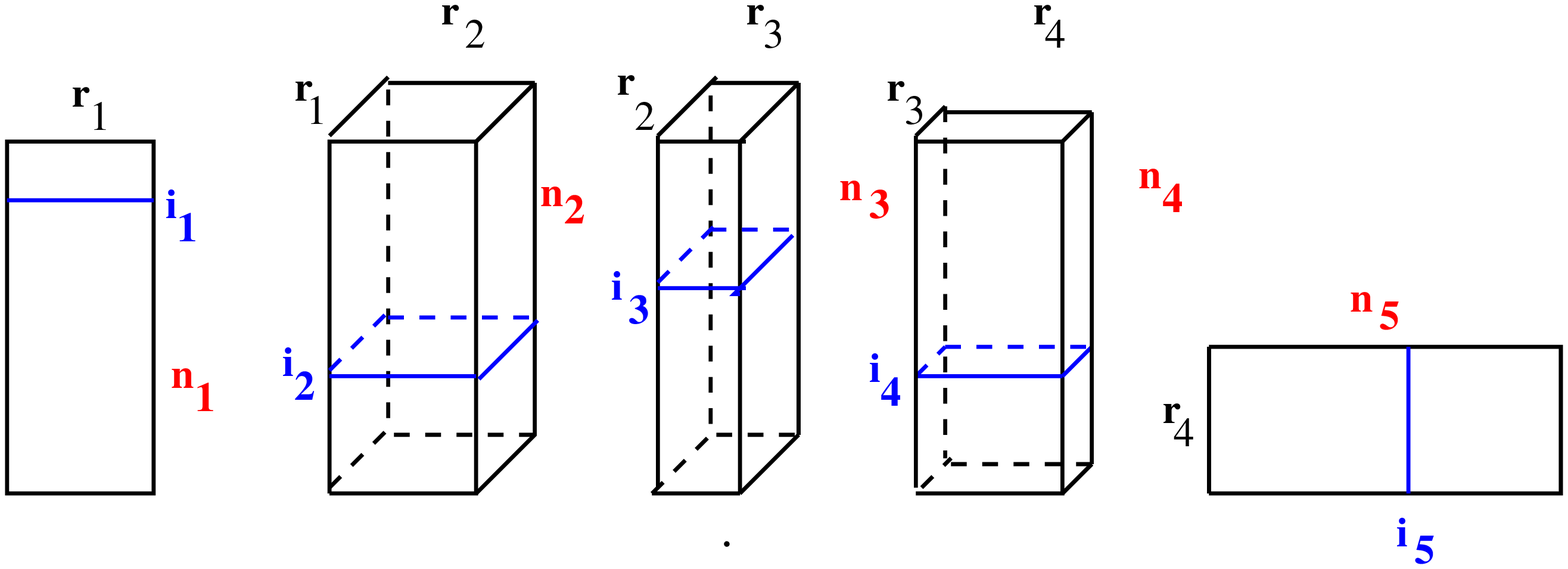}
\end{center}
\caption{Visualizing $5$th-order MPS/TT tensor.}
\label{fig:MPS_format}
\end{figure}

Figure \ref{fig:MPS_format} illustrates the TT representation of a $5$th-order tensor, where each
 particular entry  is factorized as a product of five matrices,
$a_{i_1,i_2,...,i_5} = {A}^{(1)}({i_1}) {A}^{(2)}({i_2}) \ldots {A}^{(5)}({i_5})$, where, for example,
${A}^{(2)}({i_2})\in \mathbb{R}^{r_{1} \times r_2}$.

The rank-structured tensor formats like canonical, Tucker and MPS/TT-type decompositions
also apply to matrices. For example, 
the $d$-dimensional FFT matrix over $N^{\otimes d}$ grid can be implemented on the rank-$k$ tensor  
with the linear-logarithmic cost $O(d k N\log_2 N)$,  due to the rank-$1$
factorized representation
\[
{\cal F}_{N}^{(d)}=(F_{N}^{(1)}\otimes I ...\otimes I)
(I\otimes F_{N}^{(2)}...\otimes I ) ... (I\otimes I ...
\otimes F_{N}^{(d)})
 \equiv F_{N}^{(1)}\otimes ... \otimes F_{N}^{(d)},
\]
where $ F_{N}^{(\ell)}\in \mathbb{R}^{N \times N} $ represents
the univariate FFT matrix along mode $\ell$.

\subsection{Quantics tensor approximation of functional vectors}\label{ssec:QQTform}

In the case of large mode size, the asymptotic storage cost for a class of function related
 $N$-$d$ tensors can be reduced to $O(d \log N)$ by using quantics-TT (QTT) 
tensor approximation method \cite{KhQuant:09}.
The QTT-type approximation of an $N$-vector with $N=2^L$, $L\in \mathbb{N}$, is defined as the 
tensor decomposition (approximation) in the canonical, TT or more general formats applied to a tensor 
obtained by the dyadic folding (reshaping) of the target vector to an $L$-dimensional 
$2\times ... \times 2$ data array (tensor) that is thought as an element of the $L$-dimensional 
quantized tensor space.  

In the vector case, i.e. for $d=1$, a vector 
$
{\bf x}=[x_i]\in \mathbb{R}^{{N}}, 
$
with $N=2^L$, is reshaped to its quantics (quantized) image in 
$
 \bigotimes_{j=1}^{L}\mathbb{R}^{2},
$ 
by dyadic folding, 
\[
\mathcal{F}_{2,L}: {\bf x} \to {\bf Y}
=[y_{\bf j}]\in {\bigotimes}_{j=1}^{L}\mathbb{R}^{2}, \quad {\bf j} = (j_{1},...,j_{L}), 
\;\,  j_{\nu}\in \{1,2\}, 
\]
where for fixed $i$, we have $y_{\bf j}:= x_i$, with $j_\nu=j_\nu(i)= C_{\nu -1}$ 
($\nu=1,...,L$) being defined via binary coding,
i.e. the coefficients $C_{\nu -1}\in \{0,1\} $ are found from the binary representation of $i-1$,
\[
i-1 =  C_{0} +C_{1} 2^{1} + \cdots + C_{L-1} 2^{L-1}\equiv
\sum\limits_{\nu=1}^L (j_{\nu}-1) 2^{\nu-1}. 
\]
Next figure visualizes the QTT approximation process.
\begin{figure}[tbh]
  \begin{center}
\includegraphics[width=9.0cm]{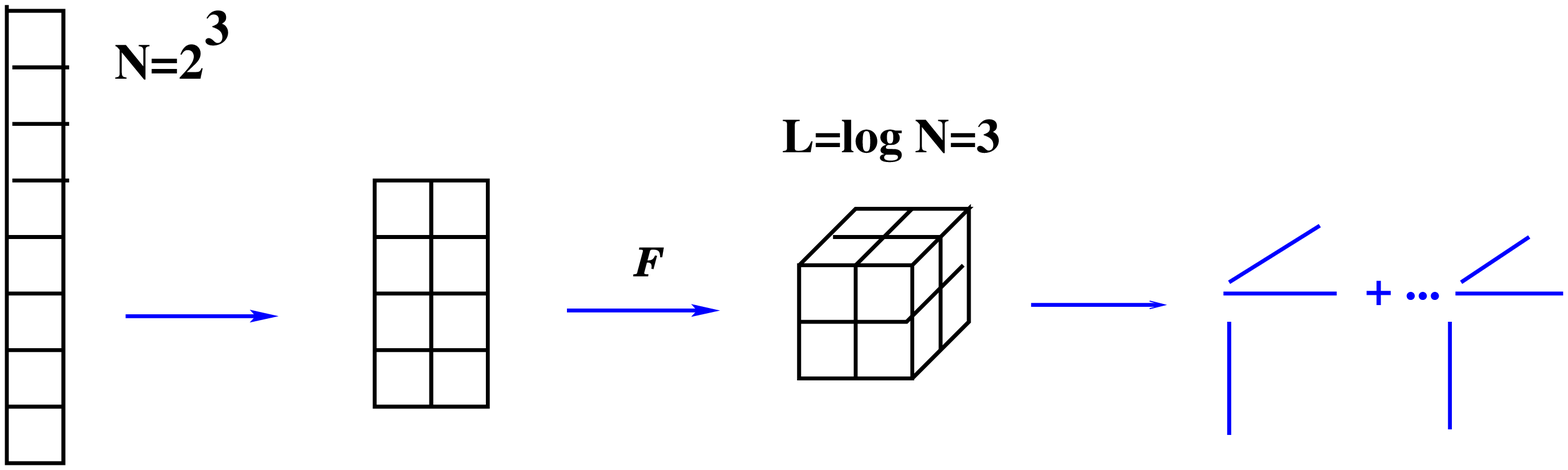}
\end{center}
\end{figure}

Suppose that the quantics image for an $N$-vector,  
i.e. an element of $L$-dimensional quantized tensor space $\bigotimes_{j=1}^{L}\mathbb{R}^{2}$
with $L=\log_2 N $,
can be represented (approximated) by the low-rank canonical or TT tensor of order $L$, 
thus introducing the QTT approximation of an $N$-vector. 
Given rank parameters $\{r_k\}$ ($k=1,...,L$), the QTT approximation of an $N$-vector  
requires a number of  representation parameters estimated by
$$
2   r^2 \log_2 N \ll N, \quad \mbox{where}\quad r_k \leq r,\quad  k=1,...,L,
$$ 
providing $\log$-volume complexity scaling in the size of initial vector, $N$.
For $d> 1$ the construction is similar \cite{KhQuant:09}.

The power of QTT approximation method is explained by the theoretical substantiation of 
the QTT approximation properties discovered in \cite{KhQuant:09} and establishing
the perfect rank-$r$ decomposition  for 
the wide class of function-related tensors obtained by sampling  continuous functions
over uniform (or properly refined) grid:
\begin{itemize}
\item $r=1$ for complex exponents, 
\item $r=2$ for trigonometric functions and plain-waves.
\item $r\leq m+1$ for polynomials of degree $m$, 
\item $r$ is a small constant for wavelet basis functions, Gaussians, etc.
\end{itemize}
all independently on the vector size $N$. 

Notice that the notion quantics (or quantized) tensor approximation (with a shorthand QTT) 
originally introduced in 2009, see \cite{KhQuant:09}, is reminiscent of  
the entity ``quantum of information'', that mimics 
the minimal possible mode size ($n=2$) of the quantized image.

Concerning the matrix case, 
it was first found in \cite{Osel-TT-LOG:09} by numerical tests 
that in some cases the dyadic reshaping of an $N \times N$ matrix with $N=2^L$ may lead to a small
TT-rank of the resultant  matrix rearranged to the tensor form.
The efficient low-rank QTT representation for a class of discrete multidimensional 
elliptic operators (matrices) and their inverse  was proven in \cite{KazKhor_1:10}.
Moreover, based on the QTT approximation, the important algebraic matrix operations
like FFT, convolution and wavelet transforms can be implemented by superfast algorithms
in $O(\log_2 N)$-complexity, see survey paper \cite{KhorSurv:14} and references therein.

\subsection{ Rank-structured tensor operations in 1D complexity}
\label{ssec:Oper_1D}

The rank-structured tensor representation provides 1D complexity of multilinear operations 
with multidimensional tensors. 
Rank-structured tensor representation  provides fast multi-linear algebra with linear 
complexity scaling in the dimension $d$. 

For given canonical tensors ${\bf A}$ and ${\bf B}$ as in (\ref{eqn:CP_form}),
with ranks $R_a$ and $R_b$, respectively, their Euclidean scalar product can be computed by
univariate operations
\begin{equation}\label{scal}
\left\langle {\bf A}, {\bf B} \right\rangle =
\sum\limits_{i=1}^{R_a} \sum\limits_{j=1}^{R_b}
c_i c_j \prod\limits_{\ell=1}^d \left\langle  {\bf a}_i^{(\ell)},  {\bf b}_j^{(\ell)} \right\rangle,
\end{equation}
at the expense $O(d n R_a R_b)$.

The Hadamard (entrywise) product of tensors ${\bf A}$, ${\bf B}$ is defined by 
${\bf Y}=[y_{i_1,...,i_d}]:={\bf A}\odot {\bf B}$, where $y_{i_1,...,i_d}= a_{i_1,...,i_d}b_{i_1,...,i_d}$.
For canonical tensors ${\bf A}$ and ${\bf B}$ given in form (\ref{eqn:CP_form}),
the Hadamard product is calculated in $O( d n R_a R_b )$ operations by 1D entrywise products of vectors, 
\begin{equation}\label{had}
{\bf A}  \odot {\bf B}  =
\sum\limits_{i=1}^{R_a} \sum\limits_{j=1}^{R_b}
c_i c_j \left(  {\bf a}_i^{(1)} \odot  {\bf b}_j^{(1)} \right) \otimes \ldots \otimes
\left(  {\bf a}_i^{(d)} \odot  {\bf b}_j^{(d)} \right).
\end{equation}
Summation of two tensors in the canonical format ${\bf C} ={\bf A}  + {\bf B}$ is performed 
by a simple concatenation of their factor matrices, 
$A^{(\ell)} =[{\bf a}_1^{(\ell)}, \ldots ,{\bf a}_{R_a}^{(\ell)}]$ and
$B^{(\ell)} =[{\bf b}_1^{(\ell)}, \ldots ,{\bf b}_{R_b}^{(\ell)}]$,
\begin{equation}\label{eqn:conc}
C^{(\ell)} =[{\bf a}_1^{(\ell)}, \ldots ,{\bf a}_{R_a}^{(\ell)},{\bf b}_1^{(\ell)}, \ldots ,
{\bf b}_{R_b}^{(\ell)}] \in \mathbb{R}^{n_\ell \times (R_a + R_b)}.
\end{equation} 
The rank of the resulting canonical tensor increases up to $R_c =R_a + R_b$.

In electronic structure calculations, the 3D convolution transform with the 
Newton kernel, $\frac{1}{\|x-y\|}$, is the most computationally expensive operation.
The tensor method to compute  convolution over large $n\times n\times n$  Cartesian grids
in $O(n\log n)$ complexity was introduced in \cite{Khor1:08}.
Given canonical tensors ${\bf A}$, ${\bf B}$,  
their convolution product  is represented by the sum of tensor products of $1D$ convolutions, 
\begin{equation}\label{conv}
 {\bf A}  \ast {\bf B} = \sum\limits_{i=1}^{R_a}
\sum\limits_{ j= 1}^{R_b}  c_i c_j \left(  {\bf a}_i^{(1)} \ast {\bf b}_j^{(1)} \right) 
\otimes \left(  {\bf a}_i^{(2)} \ast {\bf b}_j^{(2)} \right) 
 \otimes \left(  {\bf a}_i^{(3)} \ast {\bf b}_j^{(3)} \right) ,
\end{equation}
where ${\bf a}_k^{(\ell)} \ast {\bf b}_m^{(\ell)}$ denotes the univariate convolution product of $n$-vectors.
The cost of tensor convolution in both storage and time is estimated by 
$O(R_a R_b n \log n)$. The resulting algorithm considerably outperforms the 
conventional  3D FFT-based approaches of complexity $O(n^3 \log n)$, see 
numerics in \cite{KhKh3:08}.

The sequences of rank-structured 
operations on matrices and vectors normally lead to the increase of tensor ranks,
usually being multiplied or added after each operation. 
The necessary  rank reduction in the Tucker and MPS type formats can be implemented by stable algorithms 
based on the higher order SVD. In the physical community the HOSVD-type algorithms 
are known since longer as the Schmidt decomposition \cite{Vidal:03,Cirac_TC:04,Scholl:11}. 
In the case of canonical tensors the rank reduction 
can be performed by the RHOSVD algorithm based on the canonical-to-Tucker and then 
Tucker-to-canonical transforms  
described and analyzed in \cite{KhKh:06,KhKh3:08,VeBoKh_Tuck:14} 
(see \S\ref{ssec:Can2Tuck} and Appendix), which 
demonstrated the stable behavior for most of examples in the Hartree-Fock calculations
we considered so far. The stability conditions for the RHOSVD have been discussed in
\cite{KhKh3:08,VeBoKh_Tuck:14}.

\section{Tensor calculus for the Hartree-Fock equation}
\label{sec:Tensor_HF}

\subsection{Calculation of multi-dimensional  integrals}\label{ssec:Multi_Int}

Tensor-structured calculation of the multidimensional convolution integral operators
with the Newton kernel have been introduced in \cite{KhKh3:08,KhKhFl_Hart:09,VKH:exch2010},
where on the examples of the Hartree and exchange operators  in the Hartree-Fock equation,
it was  shown that calculation of the  3D and 6D convolution integrals can be
reduced to a combination of 1D Hadamard products, 1D convolutions and 1D scalar
products. 
\begin{figure}[tbh]
\centering
\includegraphics[width=2.6cm]{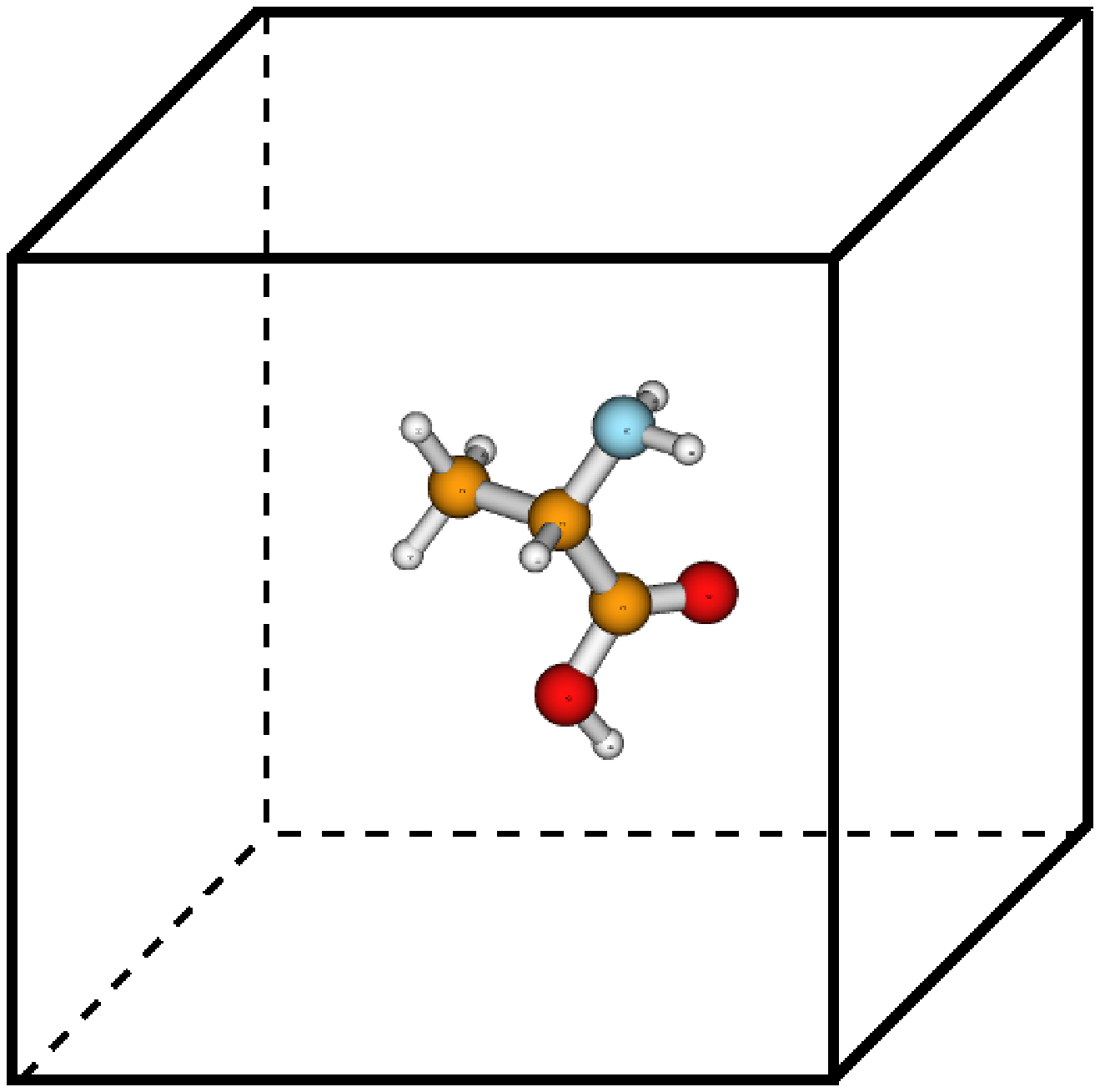}\quad\quad
\includegraphics[width=6.6cm]{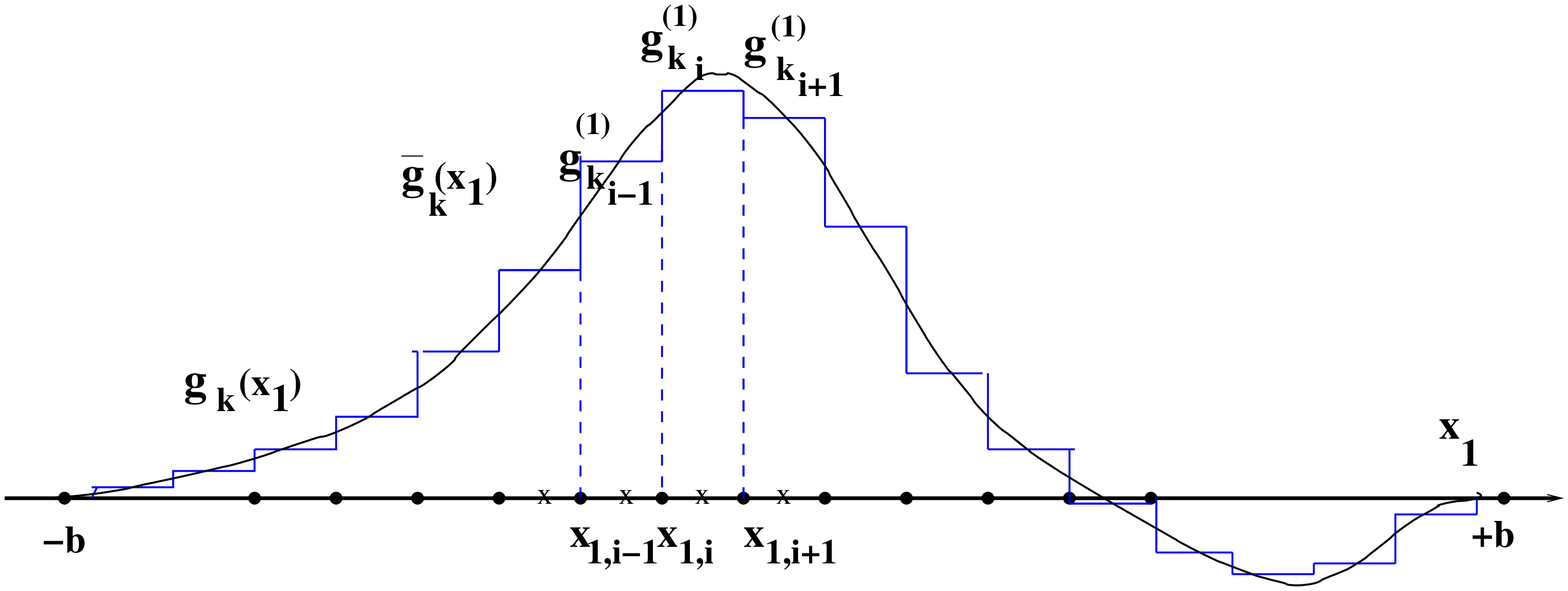}
\caption{Left: Glycine amino acid in a computational box. Right: approximation of the 
Gaussian-type basis function by a piecewise constant function.}
\label{fig:Gauss}
\end{figure}

The molecule is embedded in a certain fixed computational box $ \Omega=[-b,b]^3 \in \mathbb{R}^3$, 
as in Fig. \ref{fig:Gauss}, left. For a given discretization parameter $n \in \mathbb{N}$,  we use the 
equidistant $n\times n \times n$ tensor grid
$\omega_{{\bf 3},n}=\{x_{\bf i}\} $, ${\bf i} \in {\cal I} :=\{1,...,n\}^3 $,
with the mesh-size $h=2b/(n+1)$.
Then the Gaussian basis functions $g_k (x),\; x\in\mathbb{R}^3$,
are approximated by sampling their values at the centers of discretization intervals, as in 
Fig. \ref{fig:Gauss}, right, using one-dimensional piecewise constant basis functions
$g_k (x) \approx  \overline{g}_k (x)=\prod^3_{\ell=1} \overline{g}_k^{(\ell)} (x_{\ell})$,
$\ell=1,  2,  3$, 
yielding their rank-$1$ tensor representation,
\begin{equation}
\label{disc_GTO}
 {\bf G}_k= {\bf g}_k^{(1)} \otimes {\bf g}_k^{(2)} \otimes {\bf g}_k^{(3)}\in \mathbb{R}^{n\times n \times n},
 \quad k=1,...,N_b.
\end{equation}
Let us consider the tensor calculation of the Hartree potential
\[{V_{H}({ x}):= \int_{\mathbb{R}^3}
\frac{\rho({ y})}{\|{ x}-{ y}\|}\, d{ y}},
\]
and the corresponding Coulomb matrix,
\[
J_{km}:=\int_{\mathbb{R}^3} g_k({x}) g_m({x}) V_H({x}) d {x}, \quad
k, m =1,\ldots N_b\quad x \in \mathbb{R}^3,
\]
where the electron density, $\rho(x) = 2 \sum\limits_{a=1}^{N_{orb}}(\varphi_a )^2$,  
is represented in terms of molecular orbitals
$\varphi_a(x)=\sum\limits_{k =1}^{N_b} c_{a, k } g_k(x)$.
Given the discrete tensor representation of basis functions (\ref{disc_GTO}), the electron density is 
approximated using 1D Hadamard products of rank-$1$ tensors (instead of product of Gaussians),
\[
\rho \approx \Theta=\sum\limits^{N_{orb}}_{a=1} \sum\limits^{N_b}_{k=1} \sum\limits^{N_b}_{m=1}
c_{a, m}c_{a, k}({\bf g}_{  k}^{(1)} \odot {\bf g}_{  m}^{(1)}) \otimes ({\bf g}_{ k}^{(2)}
\odot {\bf g}_{  m}^{(2)}) 
\otimes ({\bf g}_{  k}^{(3)} \odot {\bf g}_{  m}^{(3)})
\in \mathbb{R}^{n\times n \times n}.
\]
Further, the representation of the Newton kernel $\frac{1}{\|{ x}-{ y}\|}$ by a canonical rank-$R_N$
tensor \cite{BeHaKh:08} is used (see Appendix for details),
\begin{equation} \label{eqn:NewtCan}
{\bf P}_R=
\sum\limits_{q=1}^{R_N} {\bf p}^{(1)}_q \otimes {\bf p}^{(2)}_q \otimes {\bf p}^{(3)}_q
\in \mathbb{R}^{n\times n \times n}.
\end{equation}
Since large ranks make tensor operations inefficient, the multigrid canonical-to-Tucker and 
Tucker-to-canonical algorithms  should be applied to reduce the initial rank of
$ \Theta \mapsto \Theta'$ by several orders of magnitude, from $N_b^2/2$ to $R_{\rho} \ll N_b^2/2$.
Then the 3D tensor representation of the Hartree potential is calculated by using
the 3D tensor product convolution, which is a sum of tensor products of 1D convolutions,
\[
 V_H \approx {\bf V}_H = \Theta' \ast {\bf P}_R =
\sum\limits_{ m= 1}^{R_{\rho}}
\sum\limits_{q=1}^{R_N} c_m \left(  {\bf u}^{(1)}_m \ast {\bf p}^{(1)}_q \right) \otimes
\left(  {\bf u}^{(2)}_m \ast {\bf p}^{(2)}_q \right) 
\otimes
\left(  {\bf u}^{(3)}_m \ast {\bf p}^{(3)}_q \right).
\]
The Coulomb matrix entries $J_{k m} $ are obtained by 1D scalar products
of ${\bf V}_H$ with the Galerkin basis consisting of rank-$1$ tensors,
\[
J_{km}
\approx \langle {\bf G}_k \odot {\bf G}_m, {\bf V}_H \rangle,
\quad k, m =1,\ldots N_b.
\]
The cost of 3D tensor product convolution is $O( n\log n)$ instead of $O(n^3 \log n)$ for the
standard benchmark 3D convolution using the 3D FFT. Next Table shows CPU times (sec)  for
the Matlab computation of $V_H$ for H$_2$O molecule \cite{KhKh3:08} on a SUN station using 8 Opteron
Dual-Core/2600 processors (times for 3D FFT  for $n\geq 1024$ are obtained by extrapolation).
C2T shows the time for the canonical-to-Tucker rank reduction.
\begin{table}[tbh]
\begin{center}%
\begin{tabular}
[c]{|r|r|r|r|r|r|}%
\hline
$n^3$      & $1024^3$ & $2048^3$ & $4096^3$ & $8192^3$ & $16384^3$\\
 \hline
FFT$_3$      &  $\sim 6000$ & --&--& -- & $\sim 2$ years\\
 \hline
$C\ast C$     & $8.8$ & $20.0$ & $61.0 $ & $157.5 $ & $ 299.2 $\\
\hline
{\footnotesize C2T }     & $6.9$ & $10.9$ & $20.0$ & $37.9$ & $86.0$\\
\hline
 \end{tabular}
\caption{Times (sec) for the 3D tensor product convolution vs. 3D FFT convolution.}
\end{center}
\label{tab:TConvsFFT}
\end{table}
In a similar way, the algorithm for 3D grid-based tensor-structured calculation of 
6D integrals in the exchange potential operator was introduced in \cite{VKH:exch2010},
$K_{k m}=\sum\limits^{N_{orb}}_{a=1} K_{km,a}$ with 
  \[
\label{K_ex}
K_{k m,a}:=
\int_{\mathbb{R}^3}\int_{\mathbb{R}^3} g_k( x)
\frac{\varphi_a (x)\varphi_a (y)}{|x-y|}  g_m( y) dx dy, \quad
k, m =1,\ldots N_b,
\]
the contribution from the $a$-th orbital are approximated by tensor anzats,
$$
K_{k m,a}
:= \left\langle {\bf G}_k \odot
\left[ \sum\limits_{\mu=1}^{N_b}c_{\mu a} {\bf G}_\mu \right],
\left[{\bf G}_m \odot \sum\limits_{\nu=1}^{N_b} c_{\nu a}  {\bf G}_\nu \right]
\ast {\bf P}_R  \right\rangle.
$$
Here, the tensor product convolution is first calculated for each $a$th orbital, and
then scalar products in canonical format yield the contributions to 
entries of the exchange Galerkin matrix from the $a$-th orbital.

These algorithms were employed in the first tensor-structured solver
using 3D grid-based evaluation of the Coulomb  and exchange matrices in 1D complexity  
at every step of SCF EVP iteration \cite{VeKh_Diss:10,KhKhFl_Hart:09}.
A sequence of dyadically refined 3D Cartesian grids was used 
for reducing time  in first iterations, with an $\varepsilon$ convergence criterion 
for switching to larger grids.
This is a nonstandard computational scheme avoiding
calculation of the two-electron integrals. The accuracy for small molecules like
H$_2$O and CH$_4$ was of the order of $10^{-4}$ Hartree.
Though time performance of this solver was not compatible  with the standard Hartree-Fock packages
it was the first proof of concept for the tensor numerical methods.


\subsection{3D grid-based calculation of the two-electron integrals. }\label{sssec:TEI_Int}

\noindent The basic tensor-structured Hartree-Fock solver employs the
factorized 3D grid-based calculation of the two-electron integrals tensor, 
${\bf B}=[b_{\mu\nu\kappa\lambda}]$, in the form
\begin{equation}\label{eqn:TEI_tens}
 b_{\mu\nu\kappa\lambda}=\int_{\mathbb{R}^3}\int_{\mathbb{R}^3}
\frac{g_\mu(x) g_\nu(x) g_\kappa(y)g_\lambda(y) }{\| x-y \|} dx dy 
 = \left\langle {\bf G}_\mu \odot {\bf G}_\nu ,{\bf P}_{R} \ast\left( {\bf G}_\kappa 
\odot {\bf G}_\lambda\right) \right\rangle_{n^{\otimes 3}},
\end{equation}
by using univariate tensor operations. Introduce the side matrices $G^{(\ell)}$ representing 
on the grid the full set of canonical vectors composing  the products of the Gaussian basis 
functions, $\{ {\bf G}_\mu \odot {\bf G}_\nu\}$,
\[
 G^{(\ell)} = \left[ {\bf g}^{(\ell)}_\mu \odot {\bf g}^{(\ell)}_\nu 
\right]_{1\leq \mu,\nu\leq N_b} 
\in \mathbb{R}^{n\times N_b^2}\quad \ell=1,2,3,
\]
where in most of our Hartree-Fock calculations grids of size $n^3=32\cdot 10^3$ or 
$n^3=64\cdot 10^3$ have been used.
It was found that the large matrices $G^{(\ell)}$ of size $n\times N_b^2$ 
(e.g. $N_b^2=40000$ for Alanine amino acid) can be approximated with high accuracy by low rank 
matrices with the rank parameter bounded by $R_\ell \leq N_b$ \cite{VeKhBoKhSchn:12,VeBoKhMP2:13}. 
The corresponding low-rank factorizations (``1D density fitting'') in the form
(for $\ell=1,2,3$)
\begin{equation}\label{eqn:1DDensFit}
G^{(\ell)} \cong U^{(\ell)} {V^{(\ell)}}^T, \quad 
U^{(\ell)} \in \mathbb{R}^{n\times R_\ell},
\quad V^{(\ell)}\in \mathbb{R}^{N^2_b \times R_\ell},
\end{equation}
is computed by the truncated Cholesky decomposition of the symmetric, positive definite 
$ G^{(\ell)} {G^{(\ell)}}^T $ (see \cite{VeKhBoKhSchn:12,VeBoKhMP2:13} for more details).

Based on factorization (\ref{eqn:1DDensFit}), 
the number of convolutions in (\ref{eqn:TEI_tens})
is reduced dramatically from  $N_b^2$ to $ N_b$ at most (say from $ 40000$ to $200$).
In fact, using canonical factors from the rank-$R$ canonical tensor ${\bf P}_{R}$ representing 
the Newton kernel (\ref{eqn:NewtCan}) (see (\ref{eqn:sinc_general}) in Appendix) 
we, first, precompute the set of ``convolution`` matrices  for every space 
variable $\ell, \; \ell=1,2,3$, 
\begin{equation}\label{eqn:Factor_Mk_Conv}
M_k^{(\ell)}= {U^{(\ell)}}^T ({\bf p}_k^{(\ell)}\ast_n U^{(\ell)})
\in \mathbb{R}^{R_\ell \times R_\ell}, \quad k=1,...,R,
\end{equation}
which includes the convolution products with $R_\ell \leq N_b$ column vectors of the matrix
$U^{(\ell)} \in \mathbb{R}^{n\times R_\ell}$ instead of $N_b^2/2$ convolutions in the 
initial formulation.

Given matrices $M_k^{(\ell)}$, 
then the resulting $4$-th order TEI tensor is represented in a matrix form
\[
B=mat({\bf{B}}) \cong B_{\varepsilon}:= \sum\limits^{R}_{k=1}\odot^3_{\ell=1} 
V^{(\ell)}  M_k^{(\ell)} {V^{(\ell)}}^T \in\mathbb{R}^{N_b^2\times N_b^2}.
\]
The above nonstandard factorization of the TEI matrix $B$ allows to reduce dramatically the 
computational cost of the standard Cholesky factorization schemes \cite{Higham_CHdec:90,BeMa:13} 
applied to the reduced-rank symmetric positive definite  matrix $B$.
In fact, the low-rank Cholesky decomposition of $B$ is calculated in the following sequence. First,
 the diagonal elements of $B$ are calculated as
\[
B(i,i) = \sum\limits^{R}_{k=1}\odot^3_{\ell=1} 
V^{(\ell)}(i,:)  M_k^{(\ell)} {V^{(\ell)}(:,i)}^T. 
\]
Then the selected columns of the matrix $B$, required for the rank-truncated Cholesky 
factorization scheme, are computed by the following fast tensor operations
\[
B(:,j^\ast) = \sum\limits^{R}_{k=1}\odot^3_{\ell=1} 
V^{(\ell)}  M_k^{(\ell)} {V^{(\ell)}(:,j^\ast)}^T, 
\]
leading to representation of the matrix $B$ 
in the form of rank-$R_B$ approximate decomposition \cite{VeKhBoKhSchn:12,VeBoKhMP2:13}
\[
B:= [b_{\mu\nu,\kappa \lambda}] \approx LL^T, \quad\mbox{ with } 
\quad L\in \mathbb{R}^{N_b^2 \times R_B},\; R_B\sim N_b.
\]
This algorithm is much faster than the direct Cholesky decomposition of 
the matrix $B$ with on-the-fly computation of the required column vectors.

\begin{table}[tbh]
\begin{center}%
\begin{tabular}
[c]{|r|r|r|r|r|r|}%
\hline
$n^3$    & $16384^3$ & $32768^3$ & $65536^3$ & $131072^3$  \\
\hline
mesh (bohr) &  $0.0024$ & $0.0012$ & $6 \cdot 10^{-4}$ & $3 \cdot 10^{-4}$  \\
 \hline
Had. prod.  &  $1.6$ & $3.4$ &  $12 $ & $19 $  \\
 \hline
Density fit. &  $0.5$ & $1.0$ &  $2.0 $ & $4.0$  \\
 \hline
3D Conv. time   & $69$ & $151$ & $698 $ & $496 $  \\
\hline
Chol. time  & $2.2$ & $2.2$ & $2.2$ & $2.2$  \\
\hline
 \end{tabular}
\caption{Times (sec) for 3D grid-based calculation of the directional density fitting and 
the TEI for H$_2$O molecule.}
\label{tab:DirDensFit}
\end{center}
\end{table}

Table \ref{tab:DirDensFit} represents times (sec) for 3D grid-based calculation 
of the  directional density fitting and the TEI tensor (electron repulsion integrals) 
for H$_2$O molecule in a  box $[-20,20]^3$ bohr$^3$, performed in 
Matlab on a 2-Intel Xeon Hexa-Core/2677. Time for convolution integrals in 
(\ref{eqn:Factor_Mk_Conv}) scales almost linearly in the $1D$ grid-size $n$ as expected by theory.

\subsection{Core Hamiltonian. } \label{ssec:Core}

The Galerkin representation of the 3D Laplace operator in the nonlocal Gaussian basis 
$\{g_k (x)\}_{1\leq k \leq N_b}$, $x \in \mathbb{R}^3$, 
leads to the fully populated matrix $A_g = [a_{km}]\in \mathbb{R}^{N_b \times N_b}$.
Tensor calculation of the  matrix entries $a_{km}$ for the 
discrete Laplacian $A_g$ in the separable Gaussian basis is reduced to 1D matrix 
operations \cite{VeKh_box:13} involving
the FEM Laplacian $\Delta_3$, defined on $n\times n\times n$ grid,
\[
\Delta_3  =  \Delta_1^{(1)} \otimes I^{(2)} \otimes I^{(3)} +  
I^{(1)} \otimes \Delta_1^{(2)} \otimes I^{(3)} 
 +  I^{(1)} \otimes I^{(2)} \otimes \Delta_1^{(3)},
\]
where $\Delta_1=\frac{1}{h} tridiag \{-1,2,-1\}$. Specifically, we have 
\[
{a}_{km} = \langle \Delta_3 {\bf G}_k, {\bf G}_m\rangle,
 \]
where ${\bf G}_k$ is the tensor representation 
of Gaussian basis functions using the piecewise linear finite elements.
In the case of large $n\times n\times n$ grids, this calculation can be implemented
with logarithmic cost in $n$ by using the low-rank QTT representation of 
the large matrix $\Delta_3$, see \cite{KazKhor_1:10}.

\begin{figure}[tbh]
\centering
\includegraphics[width=4.6cm]{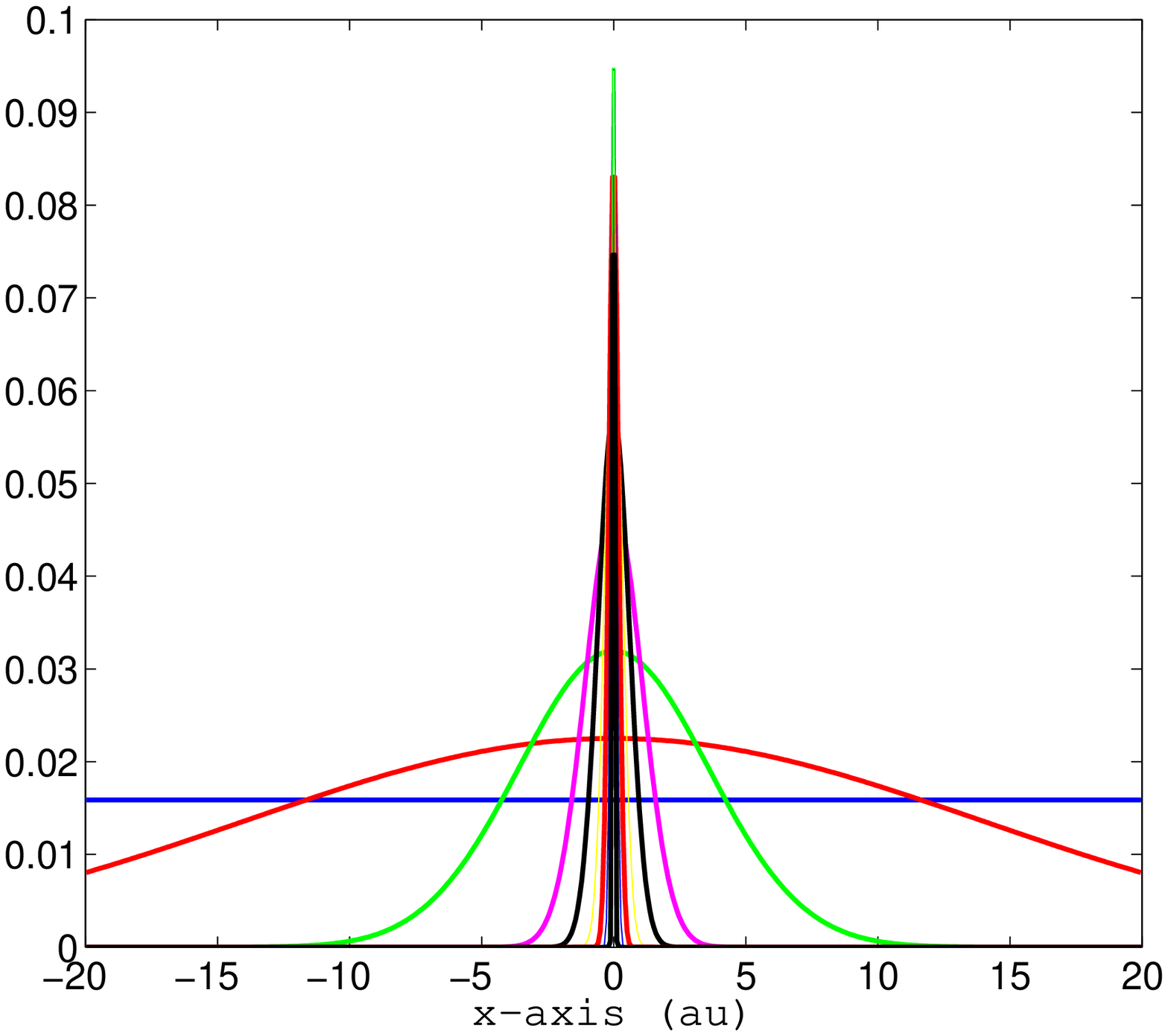}\quad
\includegraphics[width=5.2cm]{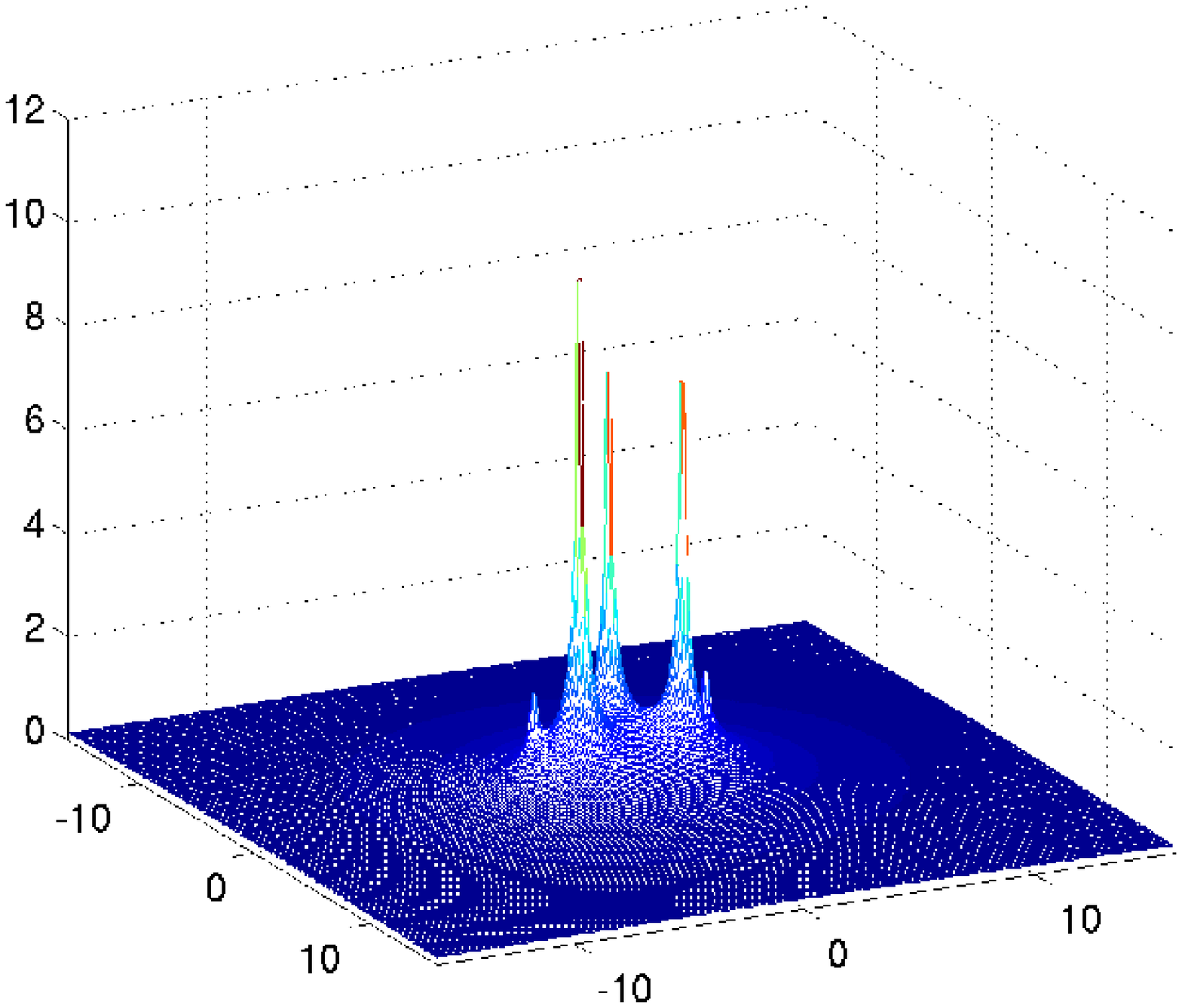} 
\caption{Left: Several vectors of the canonical representation of the Newton kernel along
one of variables. Right: calculated sum of the nuclear potentials for ethanol molecule.}
\label{fig:newt}
\end{figure}

For  tensor calculation of the nuclear potential operator  
\[
V_c(x)= - \sum_{\alpha=1}^{M}\frac{Z_\alpha}{\|{x} -a_\alpha \|},\quad
Z_\alpha >0, \;\; x,a_\alpha\in \mathbb{R}^3,
\]
we apply the rank-$1$ windowing operator, ${\cal W}_{\alpha}=
{\cal W}_{\alpha}^{(1)}\otimes {\cal W}_{\alpha}^{(2)}
\otimes {\cal W}_{\alpha}^{(3)}$, for shifting the reference Newton kernel 
$\widetilde{\bf P}_R \in \mathbb{R}^{2n\times 2n \times 2n}$ according to the coordinates
of nuclei in a molecule (see Section 5 and Appendix). 
Then the resulting nuclear potential, ${\bf P}_{c}\in \mathbb{R}^{n\times n \times n}$, 
is obtained as a direct tensor sum of shifted potentials \cite{VeKh_box:13},
\[
 {\bf P}_{c}  = \sum_{\alpha=1}^{M} Z_\alpha {\cal W}_{\alpha} \widetilde{\bf P}_R 
\]
\[
     =\sum_{\alpha=1}^{M} Z_\alpha  
\sum\limits_{q=1}^{R} {\cal W}_{\alpha}^{(1)} \widetilde{\bf p}^{(1)}_q \otimes 
{\cal W}_{\alpha}^{(2)} \widetilde{\bf p}^{(2)}_q 
\otimes {\cal W}_{\alpha}^{(3)} \widetilde{\bf p}^{(3)}_q.
\]
This leads to the following representation of the Galerkin matrix, $V_c=[\overline{v}_{km}]$, 
by tensor operations
\[
\overline{v}_{km}=  \int_{\mathbb{R}^3} V_c(x) \overline{g}_k(x) \overline{g}_m(x) dx 
\approx  \langle {\bf G}_k \odot {\bf G}_m ,   {\bf P}_{c}\rangle, 
\quad 1\leq k, m \leq N_b.
\]
Figure \ref{fig:newt}, left, shows several vectors of the canonical representation of
the Coulomb kernel along one of variables. Figure \ref{fig:newt}, right, represents the 
cross-section of the resulting nuclear potential ${\bf P}_{c}$ for C$_2$H$_5$OH molecule.

\subsection{Black-box tensor solver. }
\label{ssec:black_box}

The tensor-structured Hartree-Fock solver \cite{VeKh_box:13} based on factorized calculation of the two-electron
integrals \cite{VeKhBoKhSchn:12} includes efficient tensor implementation of the  MP2 energy 
correction \cite{VeBoKhMP2:13} scheme. Though it is yet implemented in Matlab,
its performance in time and accuracy is compatible with the standard packages based 
on analytical evaluation of the 
two-electron integrals. Due to 1D complexity of all calculations,
it enables 3D grids of the size $10^{15}$, yielding mesh size of the order of 
atomic radii, $10^{-4}$ \AA{}. That ensures high accuracy of calculations, which is 
controlled by the $\varepsilon$-ranks of tensor truncation.

The solver works in a black-box way:  input the grid-based basis-functions and 
coordinates of nuclei in a molecule and start the program. Calculation of TEI for H$_2$O
on grids $32768^3$ takes two minutes on a laptop. The time for TEI with $n^3=131072^3$
for Alanine amino acid takes approximately one hour in Matlab, including incorporated density fitting. 
 
Next examples compare the results from benchmark Molpro program 
with calculations by the tensor-based solver by using the same Gaussian basis, 
but now discretized on 3D Cartesian grid. 
For all molecules we use the ''cc-pVDZ`` Gaussian basis set.
The core Hamiltonian part in these calculations is taken from Molpro.
Note that since in the tensor solver the density fitting is included in ''blind`` calculation
of TEI, it is not easy to compare the CPU times of our calculations with those in ''ab-initio``
procedure in the standard programs, because the density fitting step there is usually considered 
as off-line pre-computing.

\begin{figure}[tbh]
\centering
\includegraphics[width=6.0cm]{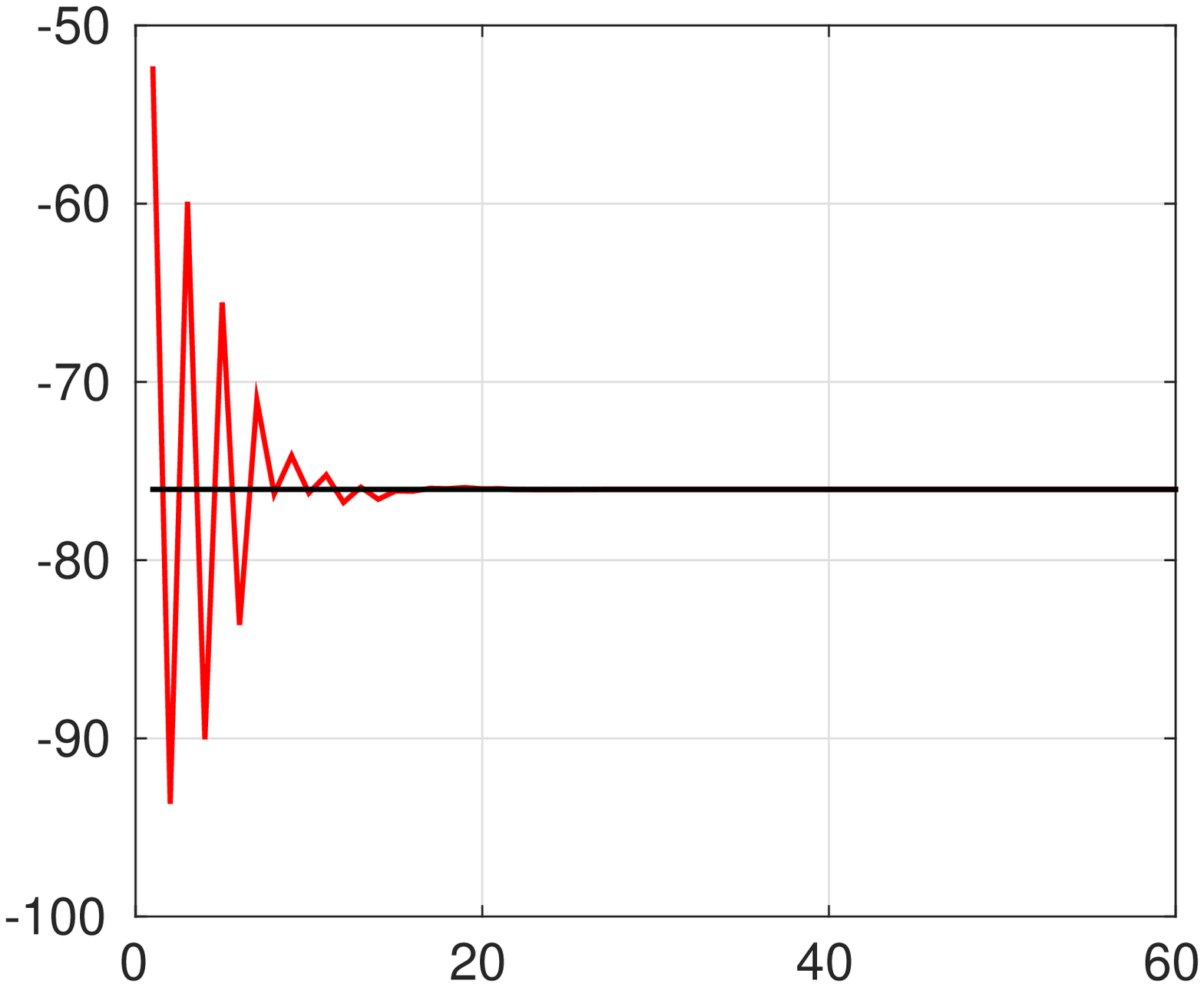}\quad
\includegraphics[width=6.0cm]{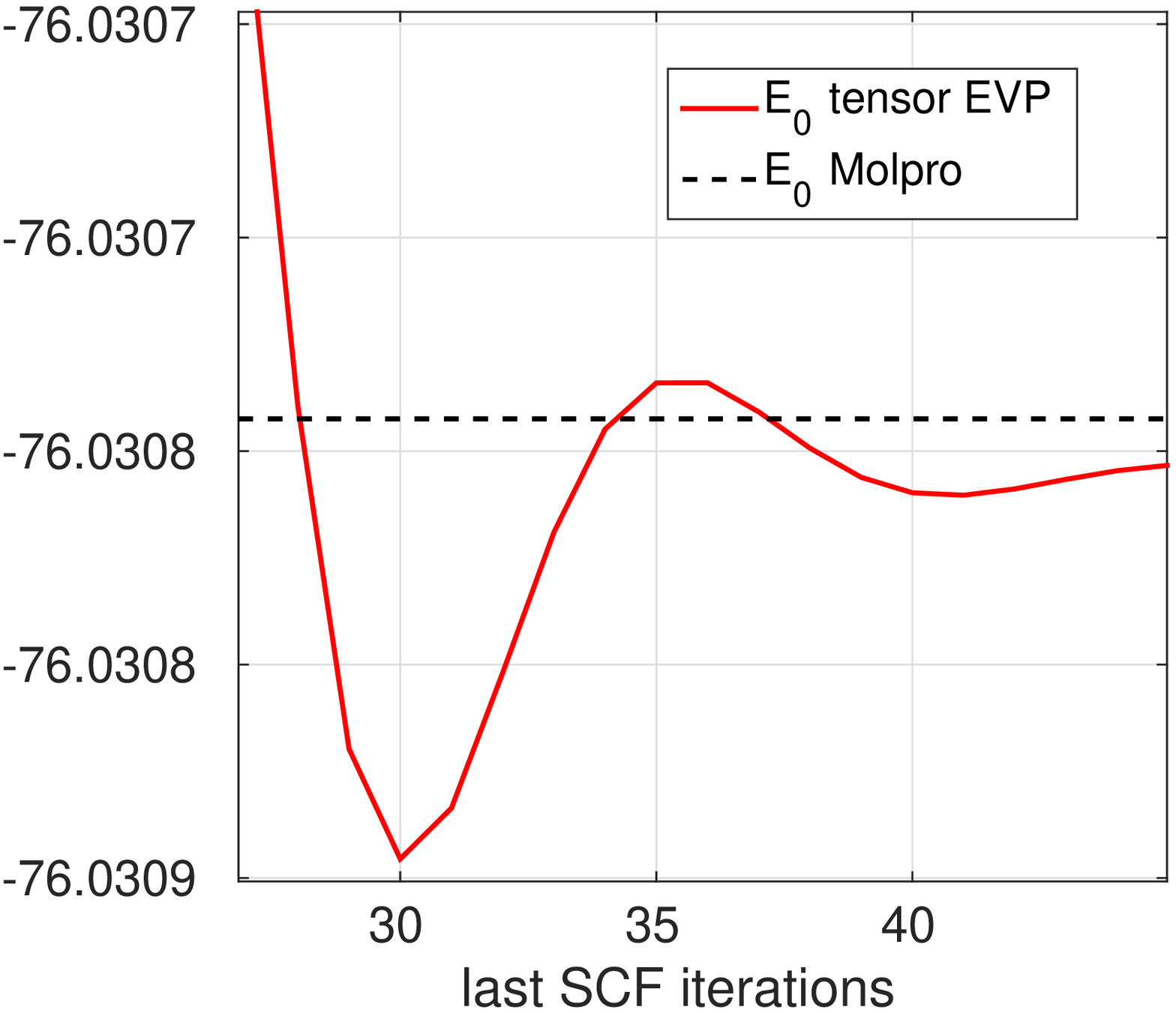} 
 \caption{Left: SCF EVP iteration for H$_2$O;
Right: convergence of the ground state energy vs. final $20$ iterations for H$_2$O.}
\label{fig:EVP_H2O}  
\end{figure}

Figure \ref{fig:EVP_H2O} shows convergence history of {\it ab initio} iterations
for H$_2$O molecule (cc-pVDZ-41), where TEI is computed on the 3D grid of size $131072^3$,
while Figure \ref{fig:EVP_H2O} (right) presents the zoom of the left graphic at the last
20 iterations. The ground state energy from Molpro is shown by the black dashed line.
 
\begin{figure}[tbh]
\centering
\includegraphics[width=6.0cm]{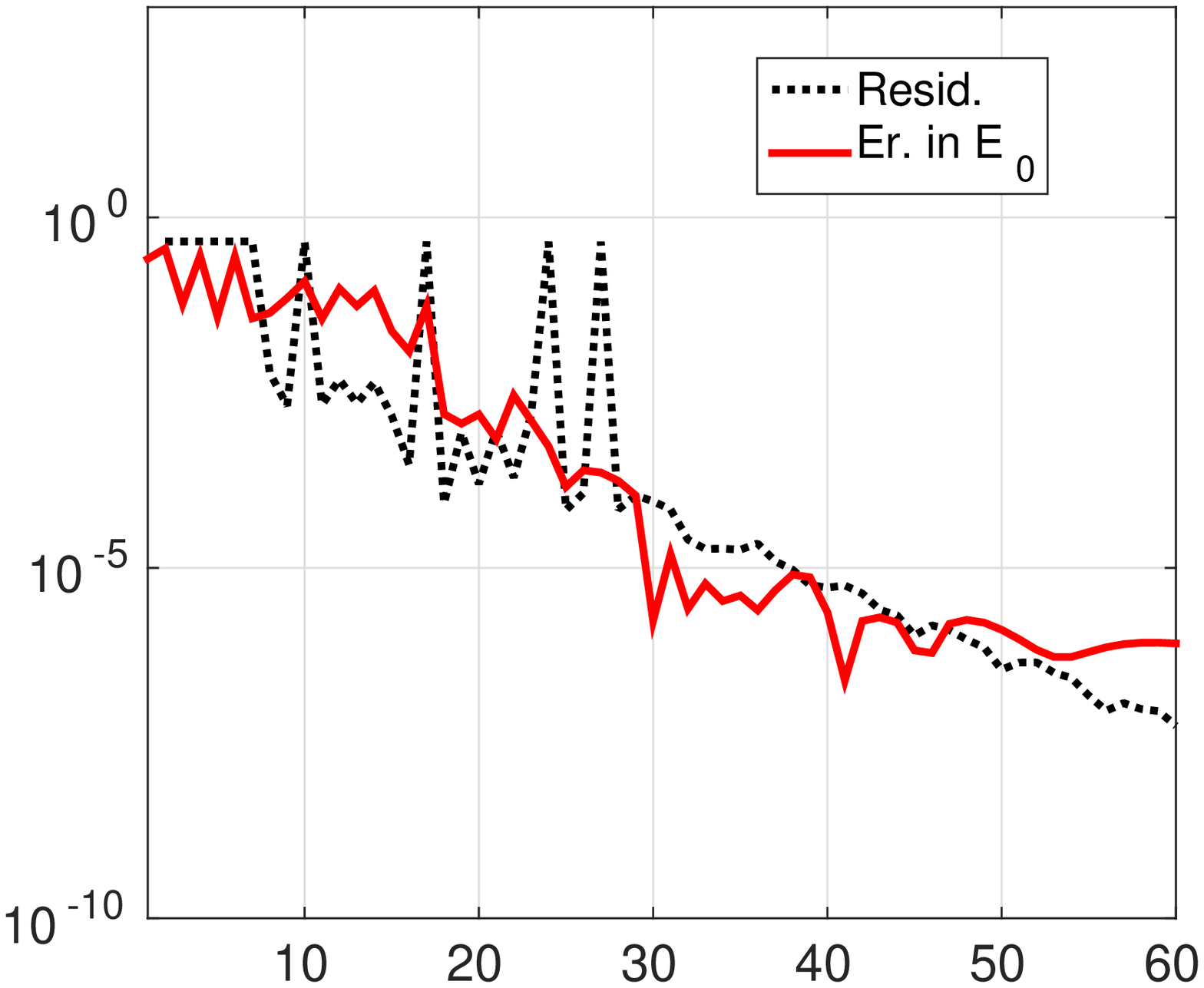}\quad
\includegraphics[width=6.0cm]{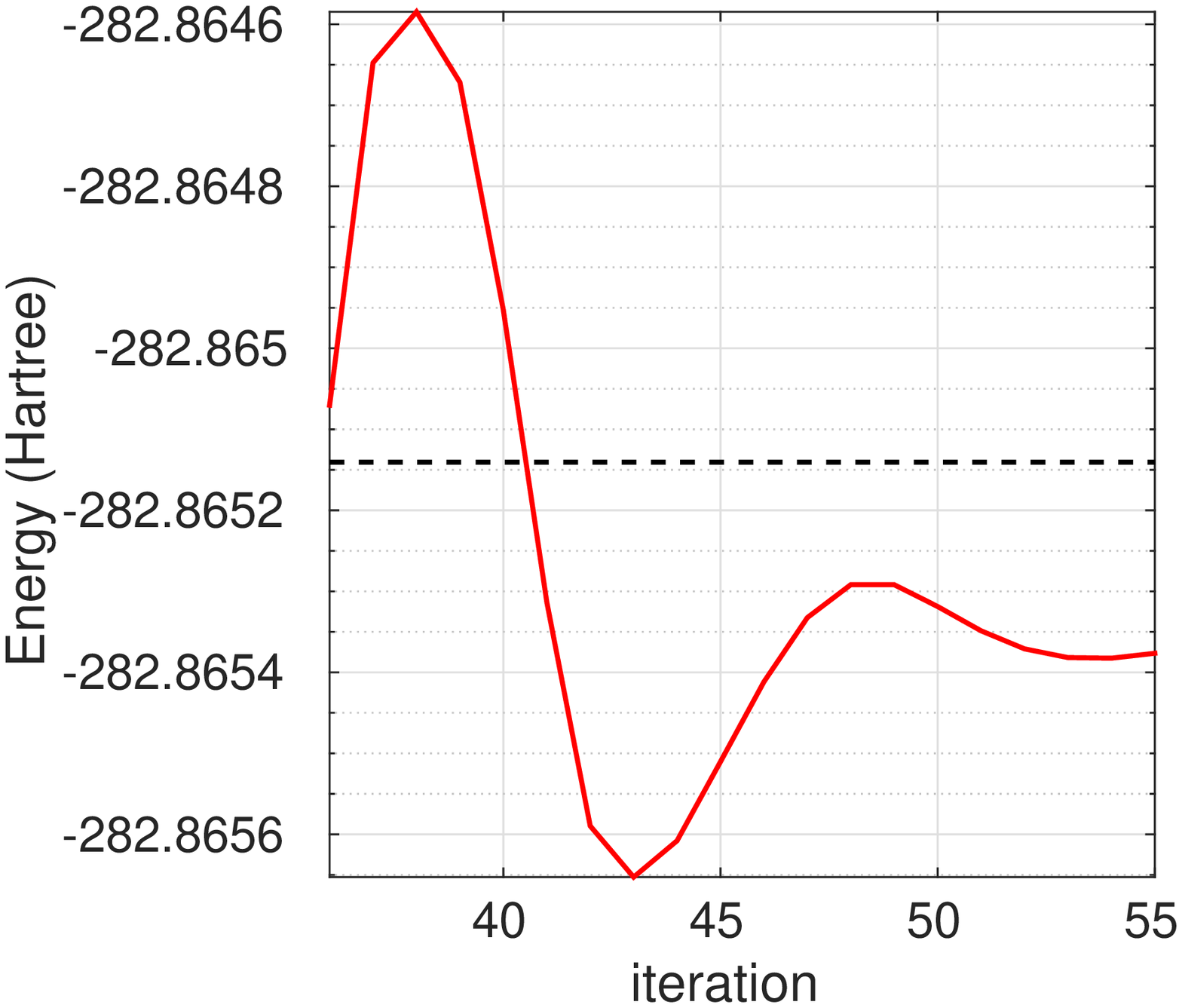} 
\caption{Left: SCF EVP iteration for Glycine (C$_2$H$_5$NO$_2$) molecule. 
Right: convergence of the ground state energy at final $20$ iterations for Glycine molecule.}
\label{fig:EVP_gly}  
\end{figure}

Figure \ref{fig:EVP_gly} (left) shows the SCF EVP iteration for Glycine amino acid 
 (cc-pVDZ-170), where dashed line indicates convergence of the residual and the red solid line 
shows convergence of the error to ground state energy from MOLPRO package with the same basis set.
Figure \ref{fig:EVP_gly} (right) illustrates the convergence of the ground state energy 
at final $20$ iterations for Glycine, dashed line is the energy from Molpro.
The second row in Table \ref{tab:FH-EVP_times} represents times (sec) for one step of 
SCF iteration in ab-initio solution of the Hartree Fock EVP for H$_2$O (cc-pVDZ-41), 
H$_2$O$_2$ (cc-pVDZ-68), and C$_2$H$_5$NO$_2$ (cc-pVDZ-170) molecules while 
the first row  shows the molecular parameters, number of orbitals and basis functions.
Next two rows show the relative difference $\Delta E_{0,g}= {|E_{0,g}-E_0|}/{|E_0|}$, 
between the grid-based ground state energy, $E_{0,g}$ and $E_0$ from Molpro with 
the same basis sets.
This numerics demonstrates that in the case of fine enough spacial $n\times n \times n$-grids
the accuracy in 7 - 8 digits (i.e. relative accuracy about $10^{-7} - 10^{-8}$)
can be achieved for moderate size molecules up to small amino acids.
\begin{table}[tbh]
\small{
\begin{center}%
\begin{tabular}
[c]{|r|r|r|r|r| }%
\hline
   & H$_2$O & H$_2$O$_2$ & C$_2$H$_5$NO$_2$   \\
\hline
 $N_{orb}$, $N_b$   & $5$; $41$ & $9$; $68$ & $20$; $170$   \\
\hline
Time  &  $0.35$ & $0.55$ & $6.0$   \\
   \hline
$\Delta E_{0,g}, (65536^3)$ & $3.0\cdot 10^{-7}$ & $8.0\cdot 10^{-8}$   & $9.1 \cdot 10^{-7}$       \\
 \hline
$\Delta E_{0,g}, (131072^3)$ & $1.4\cdot 10^{-7}$ & $3.9 \cdot 10^{-8}$  & $8.0\cdot 10^{-7}$       \\
 \hline
 \end{tabular}
\caption{Time  (sec) for one ab-initio SCF EVP iteration  and accuracy of ab-initio solution
with respect to grid-size in TEI calculations. Matlab on an Intel Xeon X5650.}
\label{tab:FH-EVP_times}
\end{center}
}
\end{table}
All calculations are performed in the computational box of size $[-20,20]^3$ bohr$^3$.
This tensor-based solver can be considered as the computational tool for trying the alternatives 
to Gaussian-type basis sets.

\section{From MP2 energy correction to excited states} \label{sec:MP2_BSE}

\subsection{MP2 correction scheme by using tensor formats}\label{ssec:MP2}

Given the set of Hartree-Fock molecular orbitals $\{C_p\}$ and the corresponding 
energies $\{\varepsilon_p\}$, $p=1,2,...,N_b$, where $\{C_i\}$ and $\{C_a\}$ 
denote the occupied and virtual orbitals, respectively.
First, one has to transform the TEI matrix  ${ B}=[b_{\mu \nu, \lambda \sigma}]$,
corresponding to the initial AO basis set,
to those represented in the molecular orbital (MO) basis, 
\begin{equation} \label{eqn:TEI_Mol}
{ V}=[v_{i a, j b}]: \quad v_{i a j b}= \sum\limits_{\mu, \nu, \lambda, \sigma =1}^{N_b} 
C_{\mu i} C_{\nu a}  C_{\lambda j} C_{\sigma b} b_{\mu \nu,  \lambda \sigma}, 
\end{equation}
where $a, b \in {\cal I}_{v}$, $i,j\in {\cal I}_{o}$, and ${\cal I}_{o}:=\{1,...,N_{orb}\}$, 
${\cal I}_{v}:=\{N_{orb}+1,...,N_{b}\}$, with $N_{orb}$ 
denoting the number of occupied orbitals.
In the following, we shall use the notation 
\[
N_{v}= N_b - N_{orb}, \quad\quad N_{ov}=N_{orb} N_{v}.
\]
The straightforward computation of the matrix $ V$ by above representation makes the 
dominating impact to the overall numerical cost of order $O(N_b^5)$. 
The method of complexity $O(N_b^4)$ based on the low-rank
tensor decomposition of the matrix $V$ was introduced in \cite{VeBoKhMP2:13}.
Indeed, it can be shown that the rank $R_B=O(N_b)$ approximation to the TEI matrix 
$
B\approx L L^T,
$ 
with the $N\times R_B$ Cholesky factor $L$, allows to introduce the low-rank representation 
of the matrix ${V}$, (see \cite{VeBoKhMP2:13} and \cite{BeKhKh_BSE:15})
\[
 V=L_V L_V^T,\quad L_V\in \mathbb{R}^{N_{ov} \times R_B},  
\]
and then reduce the asymptotic complexity of calculations to $O(N_b^4)$.

Given the tensor 
${\bf V}=[v_{i a j b}]$,  the second order MP2 perturbation to the HF energy is calculated
by
\begin{equation}\label{MP2-coor}
 E_{MP2} = - \sum\limits_{a,b \in I_{vir}} \sum\limits_{i,j \in I_{occ}} 
\frac{v_{i a j b}(2  v_{i a j b} - v_{i b j a})}{\varepsilon_a + \varepsilon_b -\varepsilon_i - \varepsilon_j}, 
\end{equation}
where the real numbers $\varepsilon_k $, $k=1,...,N_b$,
represent the HF eigenvalues. 

Introducing the so-called doubles amplitude tensor ${\bf T}$, 
\[
{\bf T}=[t_{i a j b}]: \; t_{i a j b} = 
\frac{(2  v_{i a j b} - v_{i b j a})}{\varepsilon_a + \varepsilon_b -\varepsilon_i - \varepsilon_j},
\quad a, b\in I_{vir};\; i,j\in I_{occ},
\]
the MP2 perturbation takes the form of a simple scalar product of tensors,
\[
 E_{MP2} = -  \langle {\bf V}, {\bf T} \rangle = - \langle{\bf V} \odot {\bf T}, {\bf 1} \rangle,
\]
 where ${\bf 1}$ denotes the rank-$1$ all-ones tensor. 
Introducing the low $\varepsilon$-rank reciprocal ``energy`` tensor 
\begin{equation}\label{E-tensor}
{\bf E}=[e_{abij}]:=\left[\frac{1}{\varepsilon_a + \varepsilon_b -\varepsilon_i - \varepsilon_j}\right],\quad
a, b\in I_{vir};\; i,j\in I_{occ},
\end{equation}
and the partly transposed tensor  (transposition in indices $a$ and $b$)
$$
{\bf V}'=[v'_{i a j b}]:=[v_{ibja}],
$$
allows to decompose the doubles amplitude tensor 
${\bf T}$  as follows
\begin{equation}\label{eqn:T=T1+T2}
 {\bf T}= {\bf T}^{(1)} + {\bf T}^{(2)} = 2{\bf V}\odot{\bf E} - {\bf V}'\odot{\bf E}.
\end{equation}
Notice that the denominator in (\ref{MP2-coor})
remains strongly positive if $\varepsilon_a>0$ for $a\in I_{vir}$ and 
$\varepsilon_i< 0$ for $i\in I_{occ}$. The latter condition (nonzero homo lumo gap)
allows to prove the low $\varepsilon$-rank decomposition of the tensor ${\bf E}$ 
\cite{VeKhBoKhSchn:12,VeBoKhMP2:13}.

Each term in the right-hand side in (\ref{eqn:T=T1+T2}) can be treated separately by using ranks-structured 
tensor decompositions of ${\bf V}$ and ${\bf E}$, possibly combined with various symmetries and data sparsity.
Numerical tests illustrating the tensor approach to the MP2 energy correction  
are presented in \cite{VeBoKhMP2:13}.

\subsection{Toward low-rank approximation of the Bethe-Salpeter equation for 
calculation of excited states}\label{ssec:BSE}

One of the commonly used approaches for calculation of the excited states in 
molecules and solids, along with the time-dependent DFT, is based on the solution 
of the Bethe-Salpeter equation (BSE), see for example \cite{Cas_BSE:95,RiTouSa1:13}.
The BSE approach leads to the challenging computational task on the solution of the eigenvalue 
problem for determining the excitation energies $\omega_n$, governed by
a large fully populated matrix of size $O(N^2_{ov}) \approx O(N^2_b)$,
\begin{equation} \label{eqn:BSE-GW1}
 \begin{pmatrix}
{ A}   &  { B} \\ 
{ B}^\ast  &  { A}^\ast  \\
\end{pmatrix} 
\begin{pmatrix}
 {\bf x}_n\\
{\bf y}_n\\
\end{pmatrix}
= \omega_n 
\begin{pmatrix}
{I}   &  { 0} \\ 
{ 0}  &  { -I}  \\
\end{pmatrix} 
\begin{pmatrix}
 {\bf x}_n\\ {\bf y}_n\\
\end{pmatrix},
\end{equation}
so that the computation of the entire spectrum is prohibitively expensive. 
Here the large matrix blocks of size $N_{ov}\times N_{ov}$ take a form
\[
 A=  \boldsymbol{\Delta \varepsilon} + V - \overline{W}, \quad B= {V} - \widetilde{W},
\]
where the diagonal ''energy'' matrix is defined by
\[
\boldsymbol{\Delta \varepsilon}=[\Delta \varepsilon_{ia,jb}]\in \mathbb{R}^{N_{ov}\times N_{ov}}: 
\quad \Delta \varepsilon_{ia,jb}= (\varepsilon_a - \varepsilon_i)\delta_{ij}\delta_{ab} ,
\]
while the matrices $\overline{W}=[\overline{w}_{ia,jb}]$ and $\widetilde{W}=[\widetilde{w}_{ia,jb}]$
are determined  by permutation of the so-called static screened interaction matrix
 ${W}=[{w}_{ia,jb}]$, via $\overline{w}_{ia,jb}= {w}_{ij,ab}$, and 
$[\widetilde{w}_{ia,jb}]= [{w}_{ib,aj}] $, respectively.
In turn, the forth order tensor ${\bf W}=[{w}_{iajb}]$ is constructed 
by certain linear transformations of the tensor ${\bf V}=[v_{iajb}]$, see \cite{Cas_BSE:95,RiTouSa1:13}.

A number of numerical methods for structured eigenvalue problems have been 
discussed in the literature \cite{Kressner:04,BeFa:08,DoKhSavOs_mEIG:13,LinSaadYan:14}.

The tensor approach to the solution of the partial BSE eigenvalue problem for equation (\ref{eqn:BSE-GW1})
proposed in \cite{BeKhKh_BSE:15} suggests to compute the reduced basis set by 
solving the simplified eigenvalue problem via the low-rank plus diagonal approximation
to the matrix blocks $A$ and $B$,  and then solve spectral problem for the 
subsequent Galerkin projection of the initial system (\ref{eqn:BSE-GW1}) to this reduced basis.
This procedure relies entirely on multiplication of the simplified BSE matrix with vectors.

It was demonstrated on the examples of moderate size molecules\cite{BeKhKh_BSE:15} that
a small reduced basis set, obtained by separable approximation
with the rank parameters of about several tens, allows to reveal several lowest
excitation energies and respective excited states with the accuracy about $0.1$eV - $0.02$eV.
\begin{figure}[htbp]
\centering
\includegraphics[width=6.2cm]{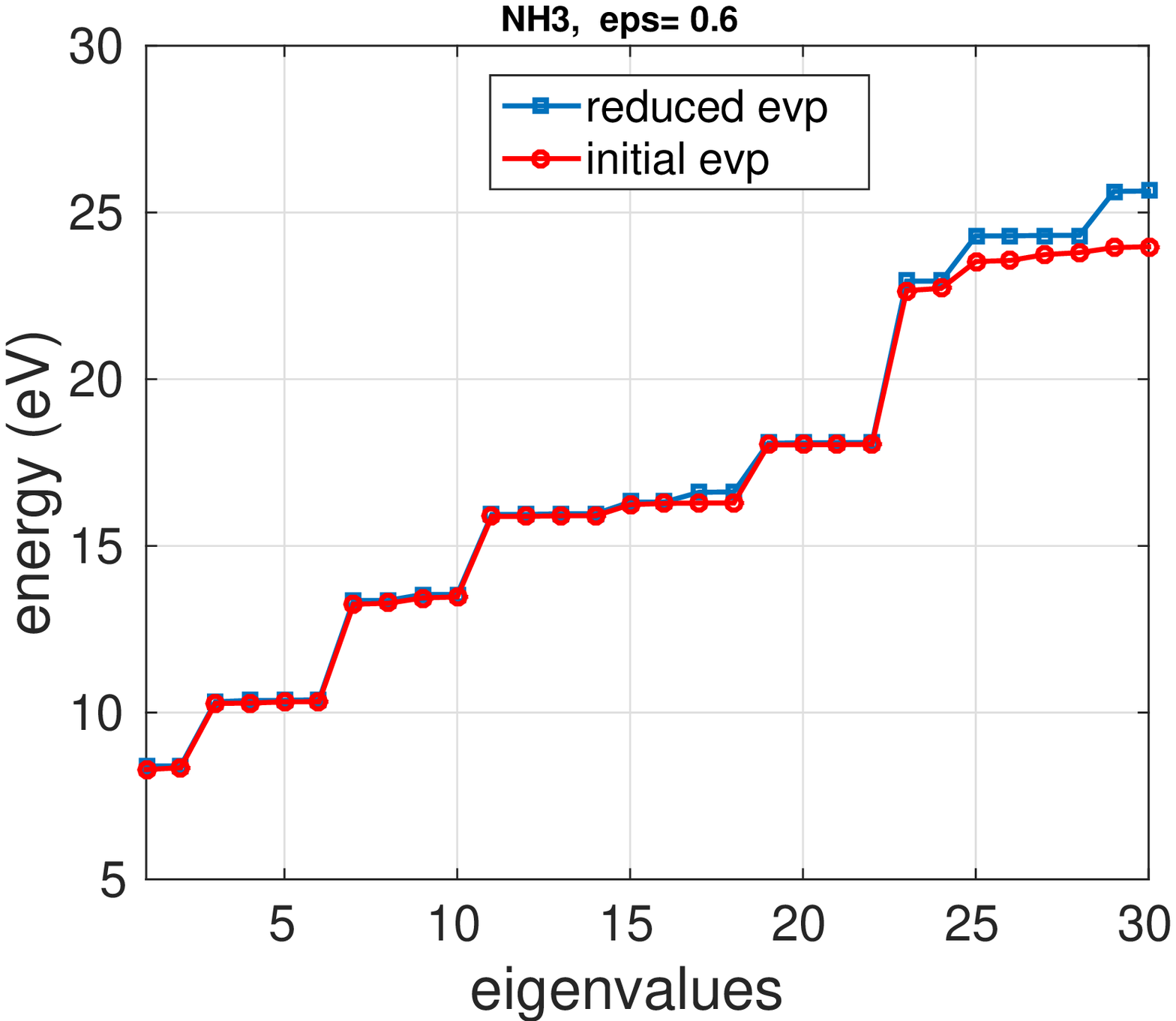} \quad
\includegraphics[width=6.2cm]{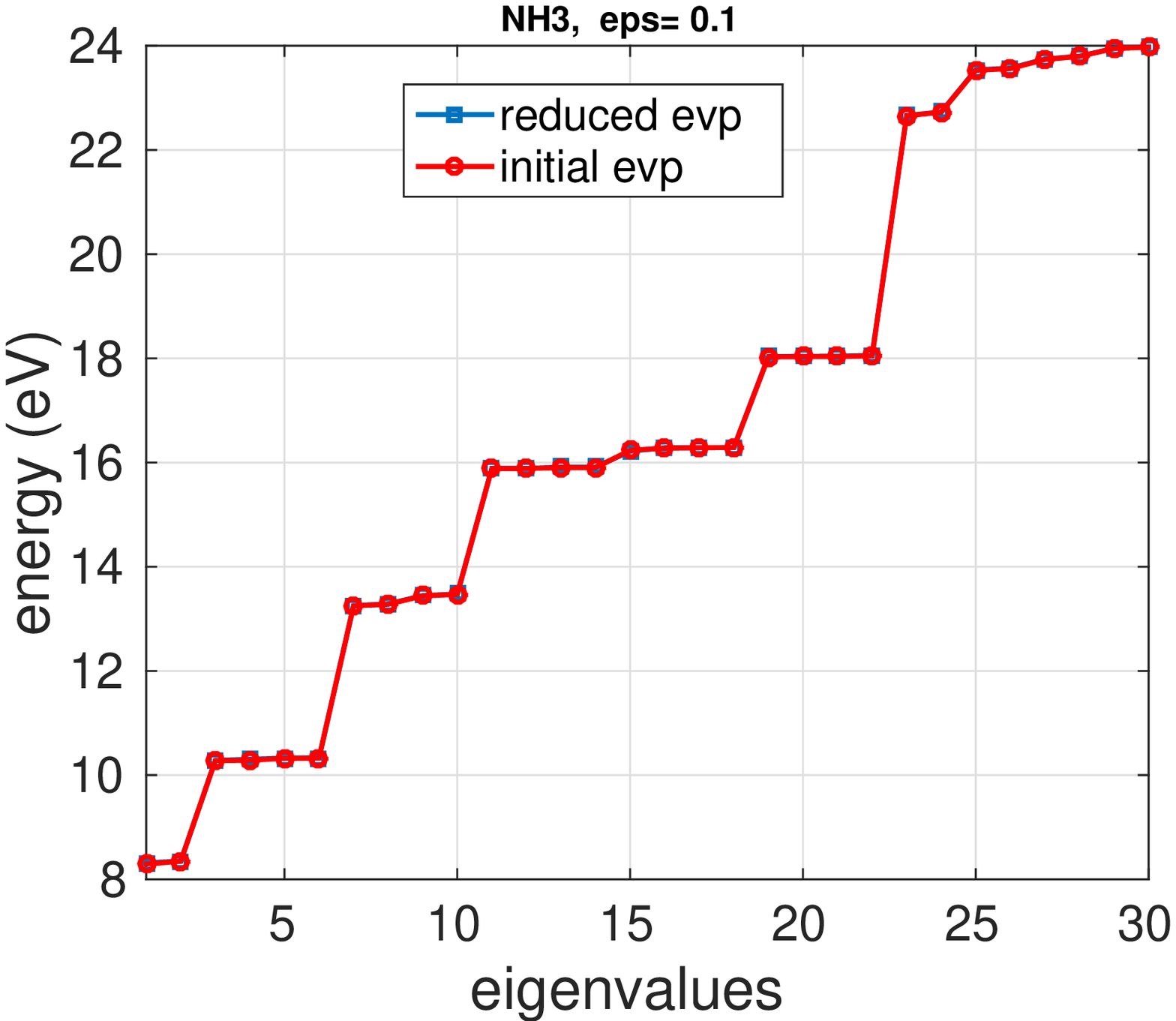}
 \caption{\small Comparison of $m_0=30$ lower eigenvalues for the reduced and exact BSE
systems for NH$_3$ molecule: $\varepsilon=0.6$, left; $\varepsilon=0.1$, right.}
\label{fig:BSE_NH3_reduced}
\end{figure}
  
Figure \ref{fig:BSE_NH3_reduced} illustrates the BSE energy spectrum of the 
NH$_3$ molecule (based on HF calculations with cc-pDVZ-48 GTO basis)
for the lowest $N_{red}=30$ eigenvalues vs. the rank truncation parameter 
$\varepsilon= 0.6$ and $0.1$, where the ranks of $V$ and the BSE matrix 
block $W$ are $4$, $5$ and $28$, $30$, respectively. 
For the choice $\varepsilon= 0.6$ and $\varepsilon=0.1$, the error in the 1st (lowest)
eigenvalue for the solution of the problem in reduced basis  is about $0.11$eV  and $0.025$eV, 
correspondingly.
The CPU time in the laptop Matlab implementation of each example is about $5$sec.

\section{Tensor approach to simulation of large crystalline clusters}
\label{sec:Cryst}

In this section, we briefly discuss 
the generalization of the tensor-based Hartree-Fock solver to the case of large lattice 
structured and periodic systems \cite{KhKh_CPC:13,KhKh_arx:14} arising 
in the numerical modeling of crystalline, metallic and polymer type compounds. 

\subsection{Fast tensor calculation of a lattice sum of interaction potentials}\label{ssec:LatticeSum}

One of the challenges in the numerical treatment of
large molecular systems is the summation of long-range potentials allocated on large 3D lattices.
The conventional Ewald summation techniques based on a separate evaluation of contributions 
from the short-  and long-range parts of the interaction potential exhibit  
$O(L^3 \log L)$ complexity scaling for a cubic $L \times L\times L$ 3D lattice.

\begin{figure}[htbp]
\centering  
\includegraphics[width=4.8cm]{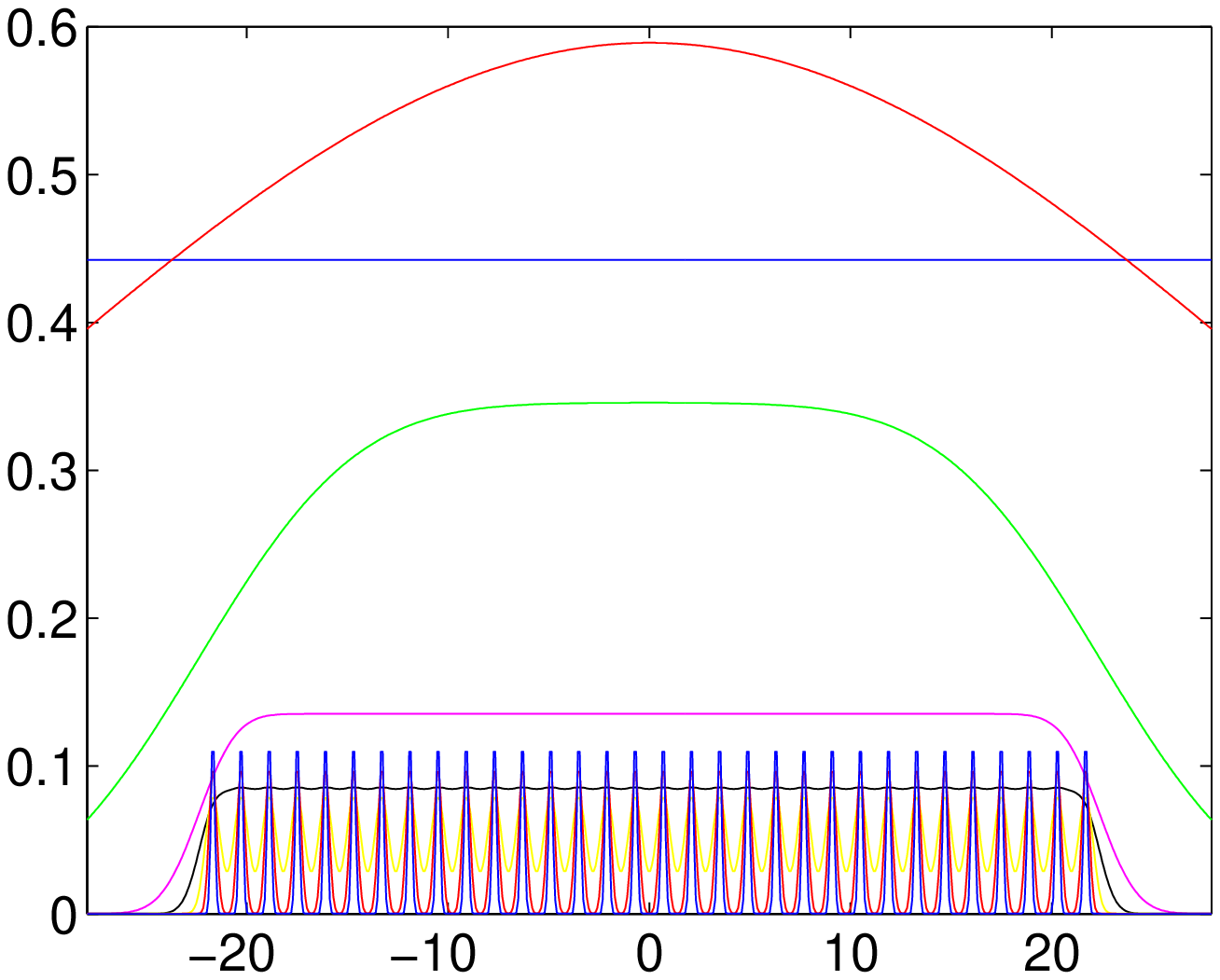}\quad
\includegraphics[width=4.8cm]{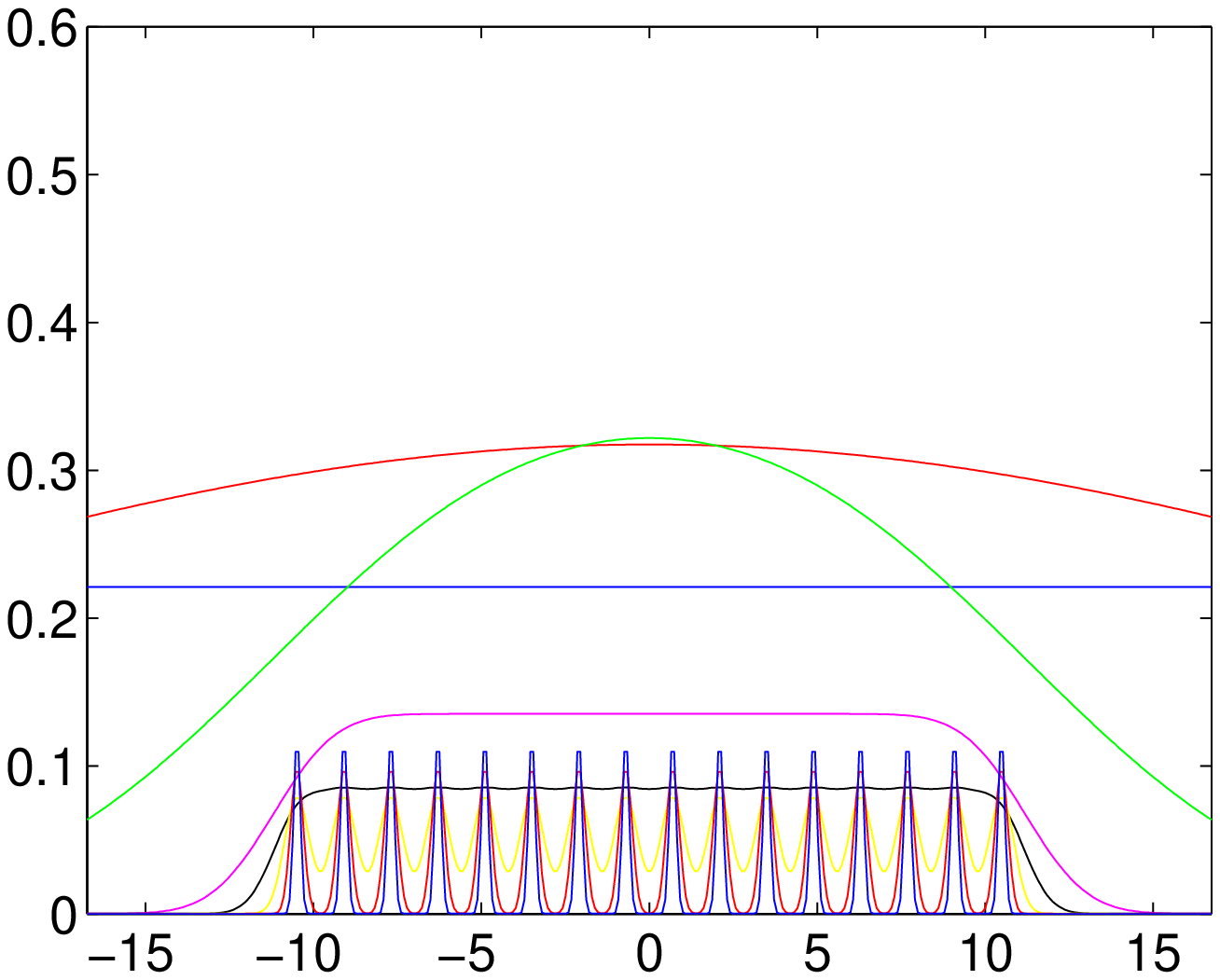}
\caption{Assembled $x$- and $y$-axis canonical vectors for a cluster of $32\times 16\times 8$ 
Hydrogen nuclei.} 
\label{fig:Assembl_can}  
\end{figure}

In the contrary, the main idea of the novel tensor summation method introduced 
in \cite{KhKh_CPC:13,VeBoKh_Tuck:14} suggests to benefit from the low-rank tensor decomposition 
of the generating kernel approximated on the fine $n\times n\times n$ representation grid
in the 3D computational box $\Omega_L$. 
This allows to completely decouple 
the 3D sum into the three independent $1$D summations, thus reducing drastically the 
numerical expenses. 
The resultant potential sum, which now requires only $O(n)$ storage and $O(n L)$  
computational demands,
is represented by a few assembled canonical/Tucker $n$-vectors of complicated shape
(see Figure \ref{fig:Assembl_can}), 
where $n$ is the univariate grid size for a cubic 3D lattice.

For ease of exposition, we consider the electrostatic potential the nuclear potential 
of a single hydrogen atom, $V_c(x)=\frac{Z}{\|x \|}$. 
Define the scaled unit cell, $\Omega_0=[0,b]^3$,   of size $b\times b \times b$ and
consider a sum of interaction potentials in a symmetric computational box   
$$
\Omega_L =B\times B \times B,\quad \mbox{with} \quad 
B= {b}[- \frac{L}{2}  ,\frac{L}{2} ], \quad L=2 L_0 \in \mathbb{N},
$$ 
consisting of a union of $L \times L \times L$ unit cells $\Omega_{\bf k}$,
obtained from $\Omega_0$ by a shift along the lattice vector $b {\bf k}$, where
${\bf k}=(k_1,k_2,k_3)\in \mathbb{Z}^3$, such that $k_\ell \in {\cal K}:={\cal K}_{-}\cup {\cal K}_{+}$ 
for $\ell=1,2,3$ with ${\cal K}_{-}:=\{-1,...,-\frac{L}{2}\}$ and 
${\cal K}_{+}:=\{0,1,...,\frac{L}{2}-1\}$. Hence, we have
$\Omega_L= \bigcup_{k_1,k_2,k_3 \in {\cal K} } \Omega_{\bf k} \in \mathbb{R}^3$.

Recall that  $b=n h$, where $h>0$ is the fine mesh size that is the same for all spatial variables,
and $n$ is the number of grid points for each variable.
We also define the accompanying domain $\widetilde{\Omega}_L$
obtained by scaling of ${\Omega}_L$ with the factor of $2$,  
 $\widetilde{\Omega}_L= 2 {\Omega}_L$, and, similar to (\ref{eqn:sinc_general}), 
introduce the respective rank-$R$ reference (master) tensor
\begin{equation}\label{eqn:refTensor}
\widetilde{\bf P}=\sum\limits_{q=1}^{R} \widetilde{\bf p}^{(1)}_{q}\otimes 
\widetilde{\bf p}^{(2)}_{q}\otimes \widetilde{\bf p}^{(3)}_{q} \in \mathbb{R}^{2n\times 2n  \times 2n},
\end{equation}
approximating $\frac{1}{\|{x} \|}$ in $\widetilde{\Omega}_L$ on a $2n\times 2n\times 2n$ 
representation grid with mesh size $h$.

Let us consider a sum of single Coulomb potentials on a $L \times L\times L$ lattice,
\begin{equation}\label{eqn:CoulombLattice}
V_{c_L}(x)=  \sum\limits_{k_1,k_2,k_3\in {\cal K}} 
\frac{Z_1}{\|{x} -a_1 (k_1,k_2,k_3)\|}, \quad x\in \Omega_L\in \mathbb{R}^3.
\end{equation}
The assembled tensor approach applies to the potentials defined on $n\times n\times n$ 
3D Cartesian grid. 
It reduces the sum over a rectangular 3D lattice, ${\bf P}_{c_L}\in \mathbb{R}^{n\times n  \times n} $,
\[
  {\bf P}_{c_L}= \sum\limits_{{\bf k} \in {\cal K}^3}  {\cal W}_{\nu({\bf k})} \widetilde{\bf P}
= \sum\limits_{k_1,k_2,k_3 \in {\cal K}} \sum\limits_{q=1}^{R}
{\cal W}_{ ({\bf k})}( \widetilde{\bf p}^{(1)}_{q} \otimes \widetilde{\bf p}^{(2)}_{q} 
\otimes \widetilde{\bf p}^{(3)}_{q}),
\]
to the summation of directional vectors for the canonical decomposition of shifted 
single Newton kernels \cite{KhKh_CPC:13}, 
\begin{equation}\label{eqn:CoulombLatTensCan}
{\bf P}_{c_L}= 
\sum\limits_{q=1}^{R}
(\sum\limits_{k_1\in {\cal K}} {\cal W}_{({k_1})}  \widetilde{\bf p}^{(1)}_{q}) \otimes 
(\sum\limits_{k_2\in {\cal K}} {\cal W}_{({k_2})} \widetilde{\bf p}^{(2)}_{q}) \otimes 
(\sum\limits_{k_3\in {\cal K}}{\cal W}_{({k_3})}  \widetilde{\bf p}^{(3)}_{q}),
\end{equation}
where ${\cal W}_{\nu({\bf k})}= {\cal W}_{{k_1}}\otimes {\cal W}_{{k_2}}\otimes {\cal W}_{{k_3}}$
is the shift-and-windowing (onto $\Omega_L$) separable transform along the ${\bf k}$-grid.
Remarkably that the rank of the resulting sum is the same as for the $R$-term 
canonical reference tensor (\ref{eqn:refTensor}) representing the single Newton kernel. 

The numerical cost and storage size are bounded by $O(R L n )$ and $O(R n)$, 
respectively, where  $n= n_0L$, and $n_0$ is the grid size in the unit cell.

\begin{figure}[htbp]
\centering  
\includegraphics[width=4.8cm]{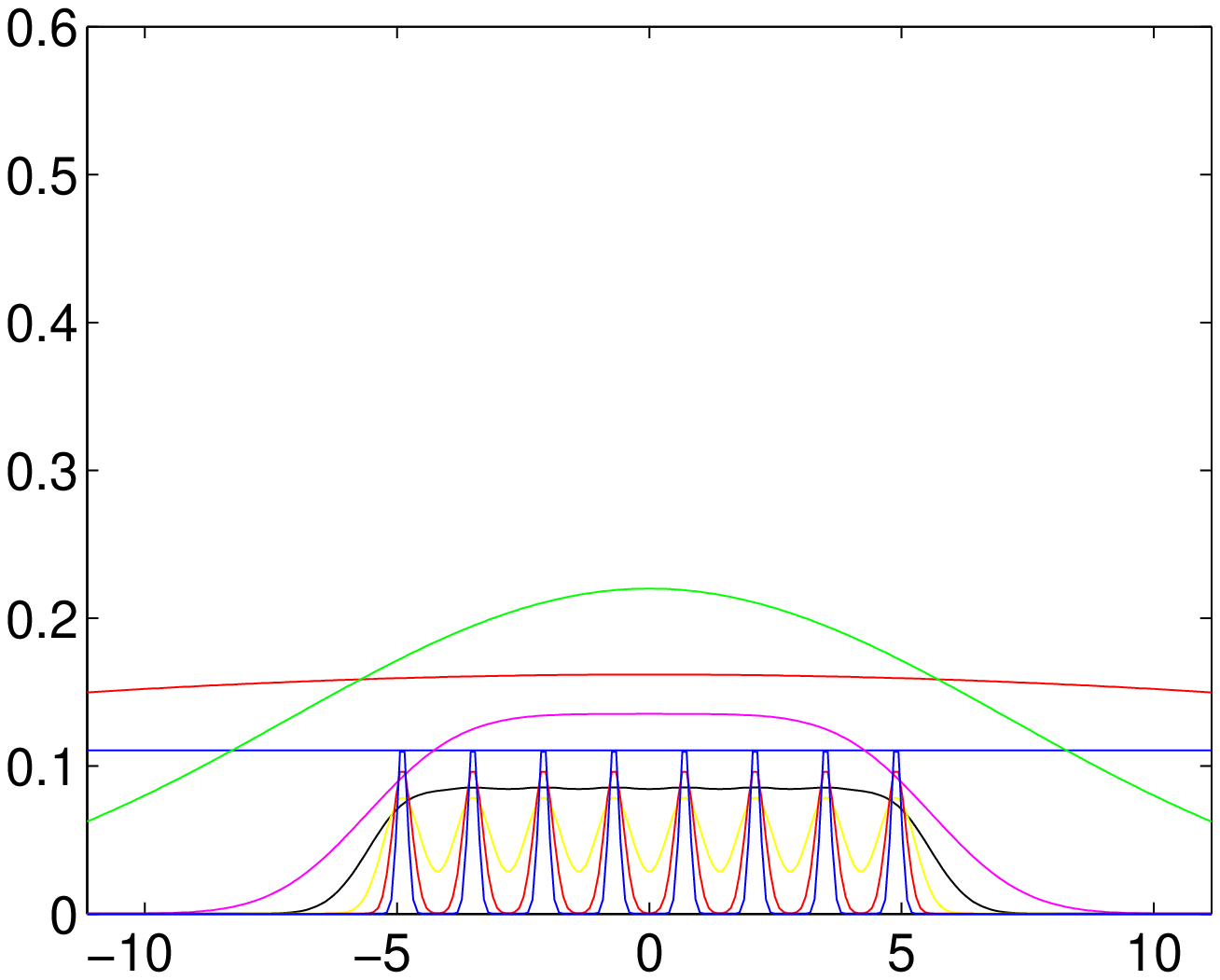}\quad
\includegraphics[width=4.8cm]{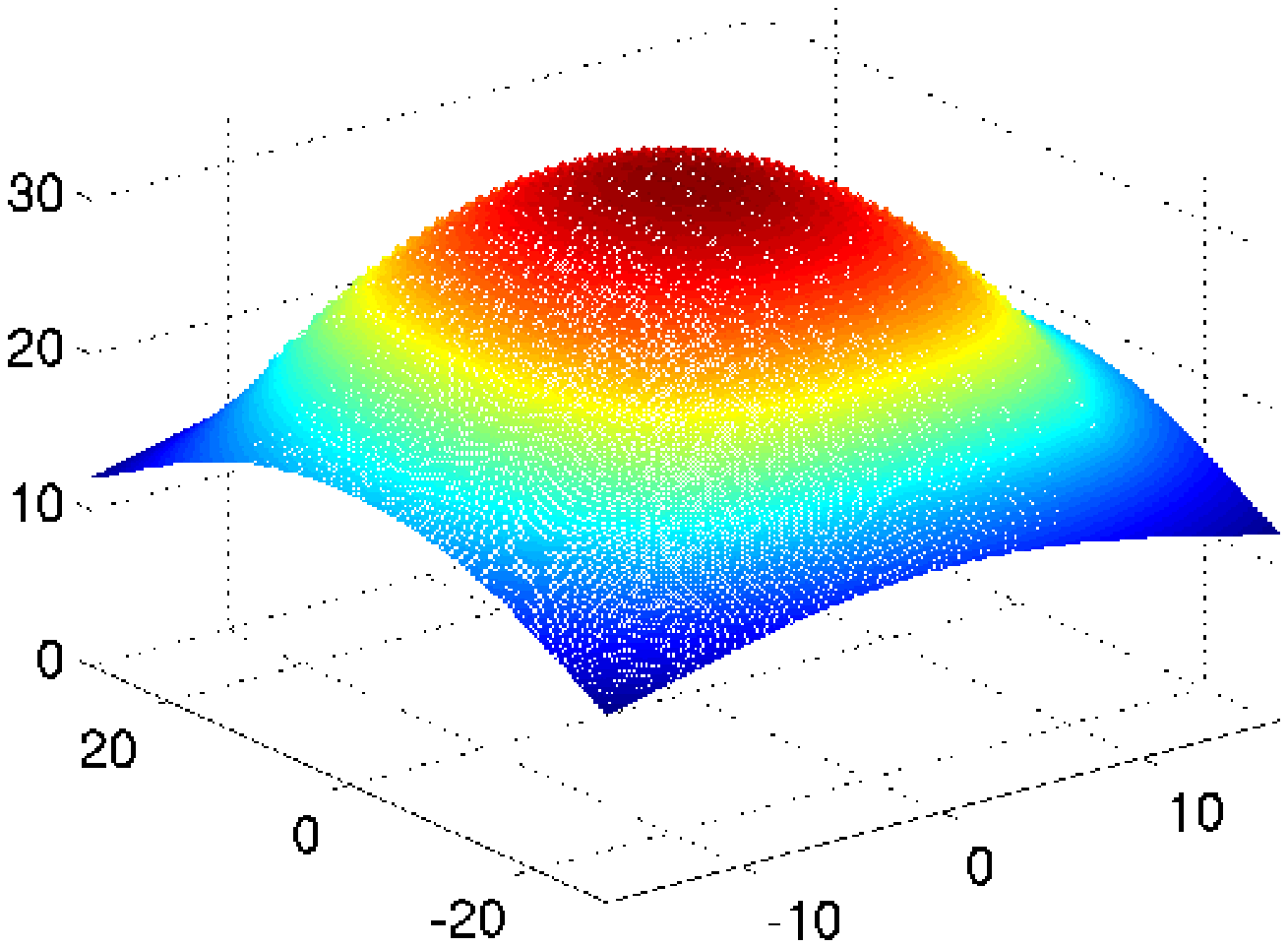}
\caption{Assembled $z$-axis canonical vectors for a cluster of $32\times 16\times 8$ 
Hydrogen nuclei. Assembled canonical sum of the Newton potentials for this cluster.} 
\label{fig:Assembl_can2}  
\end{figure}

Figures \ref{fig:Assembl_can} and \ref{fig:Assembl_can2} represent the shape of assembled canonical
vectors for the cluster of $32\times 16\times 8$ Hydrogen nuclei. 
Size of the computational box is $62 \times 32\times 22$ bohr$^3$. Here the empty interval
between the lattice and the boundary of the computational box equals to $6$ bohr.

The next table represents CPU times for the lattice summation of the Newton kernels
over $L \times L\times L$ cubic box, with very fine $n\times n\times n$ representation grid.
\begin{center}
\begin{tabular}
[c]{|r|r|r|r|r|}%
\hline
$L^3$         & $4096$  & $32768$  & $262144$ & $2097152$ \\
\hline
Time (sec.)      & $1.8$  & $0.8$  & $3.1$  &  $15.8$ \\
\hline  
3D grid size, $n^3$  & $5632^3$  & $9728^3$  & $17920^3$ & $34304^3$ \\
\hline 
 \end{tabular}
\end{center}

\begin{figure}[htbp]
\centering
\includegraphics[width=5.5cm]{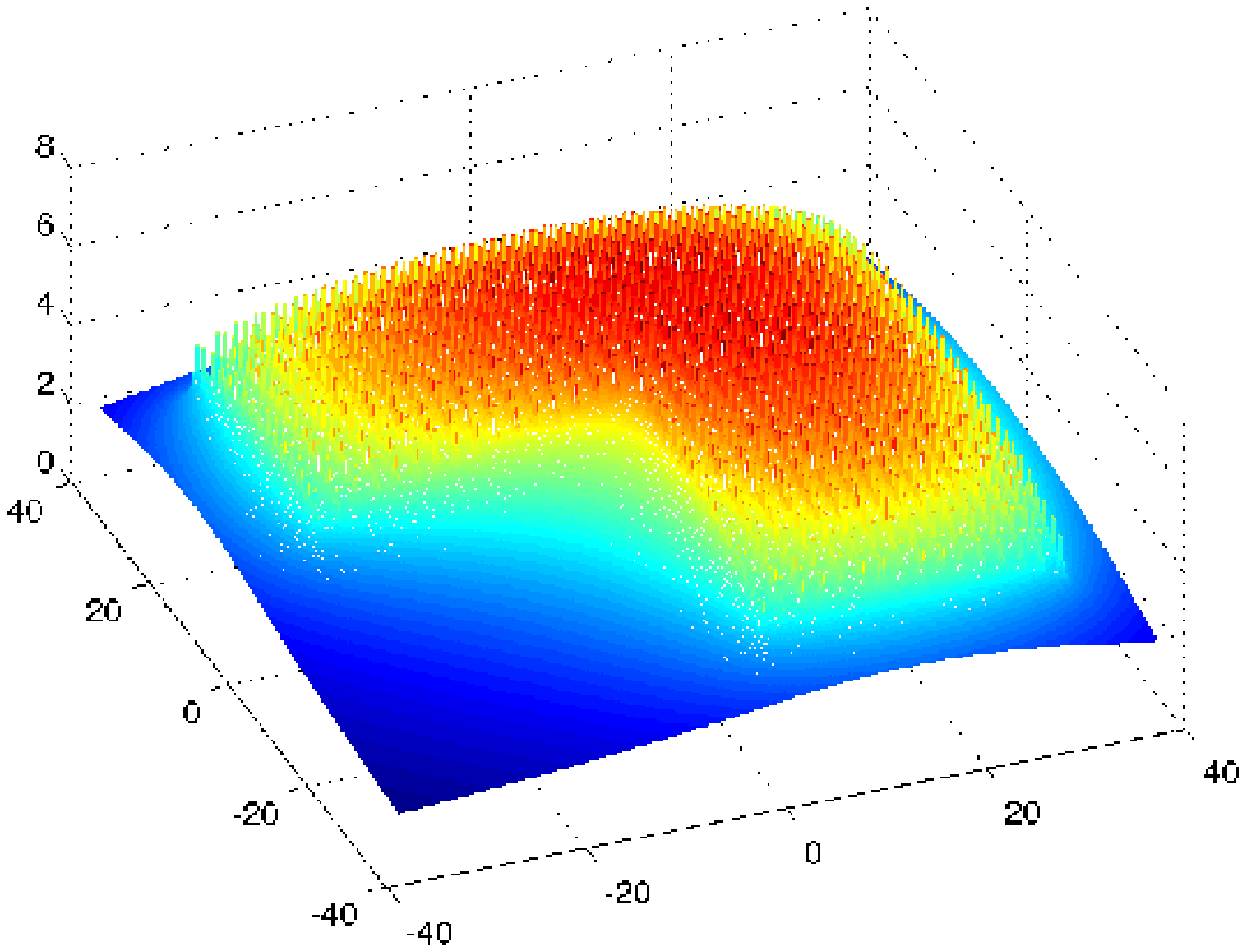}\quad
\includegraphics[width=5.5cm]{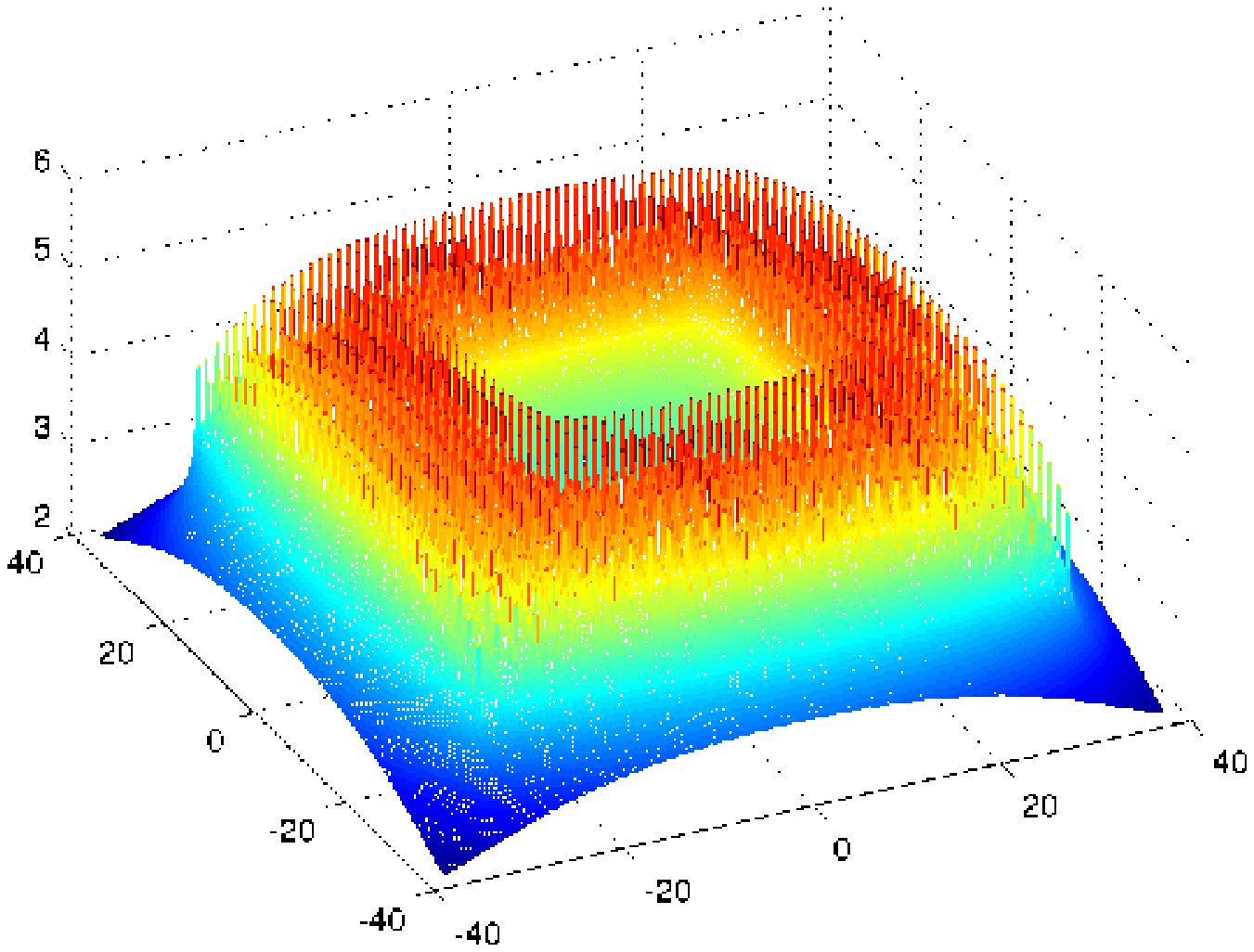}
\caption{Assembled canonical sum of the Coulomb potentials on the  $L$-shaped (left) 
and $O$-shaped (right) sub-lattices of a $32 \times 32\times 1$ lattice.}
\label{fig:Vacance_L_O}  
\end{figure}
The summation method in the canonical format was extended 
to the Tucker tensors which allows the principal generalization of this techniques 
to the case of rather complicated lattices with 
defects (Theorem 3.2, \cite{VeBoKh_Tuck:14}), so that the resulting sum takes a form 
\[ 
{\bf T}_{c_L} = 
 \sum\limits_{m_1=1}^{r_1} \sum\limits_{m_2=1}^{r_2} \sum\limits_{m_3=1}^{r_3}
 b_{m_1,m_2,m_3}
(\sum\limits_{k_1\in {\cal K} } {\cal W}_{({k_1})} \widetilde{\bf t}^{(1)}_{m_1}) 
\otimes 
(\sum\limits_{k_2\in {\cal K}} {\cal W}_{({k_2})} \widetilde{\bf t}^{(2)}_{m_2}) 
\otimes 
(\sum\limits_{k_3\in {\cal K}}{\cal W}_{({k_3})} \widetilde{\bf t}^{(3)}_{m_3}), 
\]
where $\widetilde{\bf t}^{(\ell)}_{m_\ell}$, $\ell =1,2,3$, represents the Tucker vectors of the 
rank-${\bf r}$ master tensor approximating the Newton kernel 
$\frac{1}{\|{x} \|}$ in $\widetilde{\Omega}_L$ 
on a $2n\times 2n\times 2n$ representation grid. 

In the case of defected lattices the increasing rank of the final sum of Tucker tensors 
can be reduced by the stable ALS based algorithms applicable to the Tucker tensors 
(see \cite{KhKh_arx:14}).
The particular examples of the lattice geometries suitable for our approach are
presented in Figure \ref{fig:Vacance_L_O}, see \cite{KhKh_arx:14} for more detailed discussion.

For the reduced Hartree-Fock equation, where the Fock operator is confined to the core 
Hamiltonian, the tensor-structured block-circulant representation of 
the Fock matrix was introduced\cite{KhKh_arx:14}  that allows the 
special low-rank approximation of the matrix blocks. This 
opens the way for the numerical treatment of large eigenvalue problems with
structured matrices arising in the
solution of the Hartree-Fock equation for large crystalline systems with defects.

\subsection{Interaction energy of long-range potentials on finite lattices}\label{ssec:LatticeEnergy}

Given the nuclear charges $\{{Z_{\bf k}}\}$,
centered at points ${x}_{\bf k}$,  ${\bf k}\in {\cal K}^3$, 
located on a $L \times L\times L$ lattice ${\cal L}_L=\{{x}_{\bf k}\}$ with the step-size $b$,
the interaction energy of the total electrostatic potential of these charges is 
defined by the lattice sum 
\begin{equation}\label{eqn:EnergyLatSum}
E_L = \frac{1}{2} \sum\limits_{{\bf k}, {\bf j}\in {\cal K}, {\bf k}\neq {\bf j}} 
\frac{Z_{\bf k} Z_{\bf j}}{\|{x}_{\bf j} - {x}_{\bf k}\|}, \quad 
\mbox{i.e. for} \quad \|{x}_{\bf j} - {x}_{\bf k}\|\geq b.
\end{equation}
The fast and accurate computation of electrostatic interaction energy 
(as well as the related forces and stresses) is one of the difficult tasks in 
computer modeling of macromolecular structures like finite crystal, and biological systems.

The tensor summation scheme (\ref{eqn:CoulombLatTensCan}) can be directly applied to this computational problem.
In what follows we show that (\ref{eqn:EnergyLatSum}) can be treated as a particular case of 
the previous scheme served for calculation of (\ref{eqn:CoulombLattice}) on a fine spacial grid.
For this discussion, we assume that all charges are equal, ${Z_{\bf k}}=Z$.

First, notice that the rank-$R$ reference tensor $h^{-3}\widetilde{\bf P}$ 
defined in (\ref{eqn:refTensor})
approximates with high accuracy $O(h^2)$ the (and its shifted version) Coulomb potential
$\frac{1}{\|{x} \|}$ in $\widetilde{\Omega}_L$ (for $\|{x} \|\geq b$ that is 
required for the energy expression) 
on the fine $2n\times 2n\times 2n$ representation grid with mesh size $h$.
Likewise, the tensor $h^{-3}{\bf P}_{c_L}$ approximates the potential sum $V_{c_L}(x)$ on 
the same fine representation grid including the lattice points ${x}_{\bf k}$.

We propose to evaluate the energy expression (\ref{eqn:EnergyLatSum}) 
by using tensor sums as in (\ref{eqn:CoulombLatTensCan}), but now applied to
a small sub-tensor of the rank-$R$ canonical reference tensor 
$\widetilde{\bf P}$, that is
$\widetilde{\bf P}_L:=[\widetilde{\bf P}_{|x_{\bf k}}]\in \mathbb{R}^{2L \times 2L\times 2L}$,
obtained by tracing of $\widetilde{\bf P}$ at the accompanying lattice of the double 
size $2L \times 2L\times 2L$,
i.e. ${\cal \widetilde{L}}_L=\{x_{\bf k}\}\cup \{x_{{\bf k}'}\} \in \widetilde{\Omega}_L$.
Here $\widetilde{\bf P}_{|x_{\bf k}}$ denotes the tensor entry 
corresponding to the ${\bf k}$-th lattice point designating the atomic center $x_{\bf k}$.
We are interested in the computation of the rank-$R$ tensor 
$\widehat{\bf P}_{c_L}=[{{\bf P}_{c_L}}_{|x_{\bf k}}]_{{\bf k}\in {\cal K}}\in \mathbb{R}^{L \times L\times L}$,
where ${{\bf P}_{c_L}}_{|x_{\bf k}}$ denotes the tensor entry 
corresponding to the ${\bf k}$-th lattice point on ${\cal L}_L$.
The tensor $\widehat{\bf P}_{c_L}$ can be computed at the expense $O(L^2)$ by 
\begin{equation*}\label{eqn:CoulombEnergyLatCan}
\widehat{\bf P}_{c_L}= 
\sum\limits_{q=1}^{R}
(\sum\limits_{k_1\in {\cal K}} {\cal W}_{({k_1})}  \widetilde{\bf p}^{(1)}_{L,q} \otimes 
\sum\limits_{k_2\in {\cal K}} {\cal W}_{({k_2})} \widetilde{\bf p}^{(2)}_{L,q} \otimes 
\sum\limits_{k_3\in {\cal K}}{\cal W}_{({k_3})}  \widetilde{\bf p}^{(3)}_{L,q}).
\end{equation*}
This leads to 
the representation of the energy sum (\ref{eqn:EnergyLatSum}) with accuracy $O(h^2)$ 
in a form
\[
E_L = \frac{Z^2 h^{-3}}{2} (\langle \widehat{\bf P}_{c_L}, {\bf 1}\rangle  - 
 \sum\limits_{{\bf k}\in {\cal K}} {\bf P}_{|x_{\bf k}=0} ),
\]
where the first term in brackets represents the full canonical tensor lattice sum restricted to 
the ${\bf k}$-grid composing the lattice ${\cal L}_L$, while the second term 
introduces the correction at singular points ${x}_{\bf j} - {x}_{\bf k} = 0$.
Here ${\bf 1}\in \mathbb{R}^{L \times L\times L}$  is the all-ones tensor.
By using the rank-$1$ tensor ${\bf P}_{0 L} ={\bf P}_{|x_{\bf k}=0}{\bf 1}$,
the correction term can be represented by a simple tensor operation
\[
\sum\limits_{{\bf k}\in {\cal K}} {\bf P}_{|x_{\bf k}=0} 
= \langle {\bf P}_{0 L}, {\bf 1}\rangle.
\]

Finally, the interaction energy $E_L$ allows the representation
\begin{equation}\label{eqn:EnergyLattice}
 E_L = \frac{Z^2 h^{-3}}{2} (\langle \widehat{\bf P}_{c_L}, {\bf 1}\rangle - 
\langle {\bf P}_{0 L}, {\bf 1}\rangle),
\end{equation}
that can be implemented in $O(L^2) \ll L^3 \log L$ complexity by tensor operations with
the rank-$R$ canonical tensors in $\mathbb{R}^{L \times L\times L}$. 

Table \ref{tab:EwalEnergTen} illustrates the performance of the algorithm described above.
We compare the exact value computed by (\ref{eqn:EnergyLatSum}) with the approximate 
tensor representation in (\ref{eqn:EnergyLattice}) computed on the fine representation grid
with $n= n_0 L$, $n_0=128$. We consider the lattices consisting of Hydrogen atoms with
interatomic distance $2$ bohr. The size of the largest 3D lattice with $256^3$ potentials
is of the size $256\cdot 2 +6$ bohr\footnote{Here $6$ bohr is the chosen 
dummy distance (space between the computational box boundary and the lattice).}, 
which is more than $20$ nanometers.
\begin{table}[tbh]
\begin{center}%
\begin{tabular}
[c]{|r|r|r|r|r|r|}%
\hline
$L^3$ &  Total & T$_s$ & T$_{tens.}$ & $E_L$ & err.  \\
\hline
$24^3$ & $13824$ & $37$ & $1.2$ &$3.7 \cdot 10^6$ & $2 \cdot 10^{-8}$  \\
 \hline
$32^3$ &  $32768$ & $250$ & $1.5$ &$1.5 \cdot 10^7$ & $1.5 \cdot 10^{-9}$  \\
 \hline
$48^3$ &  $110592$ & $3374$ & $2.8$ &$1.12 \cdot 10^8$ & 0  \\
 \hline
$64^3$ &  $262144$ & - & $5.7$ & $5.0  \cdot 10^8$ & -  \\
 \hline
$128^3$ &  $2097152$ & - & $13.5$ &$1.6 \cdot 10^{10}$ & -  \\
 \hline
$256^3$ &  $16777216$ & - & $68.2$ &$5.2 \cdot 10^{11}$ & -  \\
 \hline
 \end{tabular}
\caption{Comparison of times for the standard (T$_s$) and tensor-based (T$_{tens.}$)
calculation of the interaction energy for the lattice electrostatic potentials.
Matlab on an Intel Xeon X5650.}
\label{tab:EwalEnergTen}
\end{center}
\end{table}

The presented approach for fast calculation of the interaction energy can extended to the 
case of non-uniform rectangular lattices and,
under certain assumptions, to the case of non-equal nuclear charges $Z_{\bf k}$.
Moreover, it applies to many other types of spherically symmetric interaction potentials, for example, 
to shielded Coulomb interaction or van der Waals attraction sums corresponding to the distance function
$\| x\|^{-2}$ and $\| x\|^{-6}$, respectively.

\section{Conclusions}\label{sec:Conclusions}

The goal of this paper is to attract  interest of the specialists in
computational quantum chemistry to recent results and open questions
of the grid-based tensor approach in electronic structure calculations.
Here we focus mostly on the description of main mathematical and algorithmic
aspects of the tensor decomposition schemes and demonstrate their 
benefits in some applications. 

The scope of applications which can be regarded as consistent, 
ranges from the Hartree-Fock energy for moderate size molecules,
including ``blind'' calculation of TEI tensor with incorporated algebraic density fitting,
to calculation of the excited states for molecules, and up to a unique 
superfast method for calculating the lattice potential sums and the interaction 
energy of long range potentials on a lattice in a finite volume.
Tensor approach allows to treat above problems using moderate computational
facilities. All numerics given in the paper presents implementations in Matlab.

The described numerical tools are not restricted to the applications presented here, 
but can be applied to various hard computational problems in (post) Hartree-Fock 
calculations related to accurate evaluation of multidimensional integrals, 
and efficient storage and manipulations with large multivariate data arrays.


The presented method for summation of long-range potentials with sub-linear computational
cost does not have analogues in what is used so far in computational quantum chemistry 
and can have a good future. For example, it can be useful in modeling of large
finite molecular structures like nanostructures or quantum dots, where the periodic
approach may be inconsistent. Calculating a sum of several millions of
lattice potentials on fine 3D grids takes only several seconds in Matlab on a 
laptop \cite{KhKh_CPC:13}.
This method gives also the unique possibility to present the summed 3D lattice potential 
in the whole computational region with very high accuracy. Integration and differentiation of 
this 3D potential can be easily performed on representation grid due to 1D computational costs.
 
Tensor approach is now being evolved for modeling the electronic structure
of finite crystalline-type molecular systems \cite{KhKh_arx:14}. We hope that tensor numerical 
methods will have a good future in solving challenging multidimensional problems 
of computational quantum chemistry. 

\section{Appendix}\label{sec:Appendix}

{\bf 1. Canonical-to-Tucker transform}. 

The Canonical-to-Tucker tensor transform combined with the Tucker-to-Canonical scheme
introduced in \cite{KhKh3:08} usually applies 
for the rank reduction of the function related canonical tensors with the large initial rank.
Here we sketch Algorithm {\it Canonical-to-Tucker} which includes the following basic steps:

{\it Input data}: Side matrices $U^{(\ell)}=[{\bf u}_1^{(\ell)}\ldots {\bf u}_R^{(\ell)}] 
\in \mathbb{R}^{n_\ell \times R}$, $\ell=1, 2,3$,
composed of vectors ${\bf u}_k^{(\ell)}\in \mathbb{R}^{n_\ell}$, $k=1,\ldots R$, see (\ref{eqn:CP_form});
maximal Tucker-rank parameter ${\bf r}$; maximal  number of the alternating least square 
(ALS) iterations $m_{max}$ (usually a small number).

  {\bf (A)}  Compute the singular value decomposition (SVD) of side matrices:
\[
 U^{(\ell)} = Z^{(\ell)} S^{(\ell)} V^{(\ell)}, \quad \ell=1, 2, 3.
\]
Discard the singular vectors in $ Z^{(\ell)}$ and the respective singular values 
up to given rank threshold, yielding the small orthogonal matrices 
$Z^{(\ell)}_{r_\ell} \in \mathbb{R}^{n_\ell \times r_\ell}$, $\ell=1, 2, 3$.
  
 {\bf (B)} Project side matrices $U^{(\ell)}$ onto the orthogonal basis set defined by $Z^{(\ell)}_{r_\ell}$:
\begin{equation}\label{eq:alg1}
U^{(\ell)} \mapsto  \widetilde{U}^{(\ell)} = (Z^{(\ell)}_{r_\ell})^T U^{(\ell)}, \quad
 \widetilde{U}^{(\ell)}\in \mathbb{R}^{r_\ell \times R}, \quad \ell=1, 2,3.
\end{equation}

  {\bf (C)} (Find dominating subspaces). Implement the following ALS iteration $m_{max}$ times at most.
For $\ell= 1, 2,3$ implement the following ALS iteration $m_{max}$ times at most.

{\bf (D)} Start ALS iteration for  $\ell= 1, 2,3$:

$\vartriangleright$  For $\ell=1$ : construct partially projected image of the full tensor,
\begin{equation}\label{eq:alg2}
{\bf U} \mapsto
  \widetilde{\bf U}_1 = {\sum}_{k =1}^{R} c_k  {\bf u}_k^{(1)}  
\otimes \widetilde{\bf u}_k^{(2)} \otimes \widetilde{\bf u}_k^{(3)}, \quad  c_k \in \mathbb{R}.
\end{equation}
Here ${\bf u}_k^{(1)} \in \mathbb{R}^{n_1}$ are in physical space for mode $\ell=1$,
while  $\widetilde{\bf u}_k^{(2)} \in \mathbb{R}^{r_2}$ and $\widetilde{\bf u}_k^{(3)} \in \mathbb{R}^{r_3}$,
the column vectors of $\widetilde{U}^{(2)}$ and $\widetilde{U}^{(3)}$, respectively, belong to the
coefficients space by means of projection.

$\vartriangleright$ Reshape the tensor $\widetilde{\bf U}_1\in \mathbb{R}^{n_1 \times r_2 \times r_3}$ into
a matrix $M_{U_1} \in \mathbb{R}^{n_1 \times (r_2  r_3)}$, representing the span of the optimized subset 
of mode-$1$ columns in partially projected tensor $\widetilde{\bf U}_1$.
 Compute the SVD of the matrix $M_{U_1}$:
\[
 M_{U_1} = Z^{(1)} S^{(1)} V^{(1)},
\]
and truncate the set of singular vectors in 
$Z^{(1)} \mapsto \widetilde{Z}^{(1)}\in \mathbb{R}^{n_1 \times r_1}$,
according to the restriction on the mode-$1$ Tucker rank, $r_1$.

$\vartriangleright$ Update the current approximation to the mode-$1$ dominating subspace 
$Z^{(1)}_{r_1} \mapsto \widetilde{Z}^{(1)}$.

$\vartriangleright$ Implement the single loop of ALS iteration for mode $\ell=2$ and for $\ell=3$.

$\vartriangleright$ End of the single ALS iteration step.

$\vartriangleright$ Repeat the complete ALS iteration $m_{max}$ times to obtain the optimized 
Tucker orthogonal side matrices $\widetilde{Z}^{(1)}$, $\widetilde{Z}^{(2)}$, $\widetilde{Z}^{(3)}$,
and final projected image $\widetilde{\bf U}_3$.
 
{\bf (E)}  Project the final iterated tensor $\widetilde{\bf U}_3$ in (\ref{eq:alg2}) 
using the resultant basis set in $\widetilde{Z}^{(3)}$ to obtain the core tensor, 
$\boldsymbol{\beta} \in \mathbb{R}^{r_1\times r_2 \times r_3} $. 
  
{\it Output data}: The Tucker core tensor  $\boldsymbol{\beta}$ and the Tucker orthogonal
side matrices $\widetilde{Z}^{(\ell)}$, $\ell=1, 2,3$.

\vspace{2mm}
The {\it Canonical-to-Tucker} algorithm can be easily modified to the $\varepsilon$-truncation stopping criteria.
Notice that in the case of equal Tucker ranks, $r_{\ell}=r$, maximal canonical rank of the core 
tensor $\boldsymbol{\beta}$ does not exceed $r^2$, see \cite{KhKh3:08}, which
completes the Tucker-to-Canonical part of the total 
algorithm. (Further optimization of the canonical rank in the small-size core tensor 
$\boldsymbol{\beta}$ can be implemented by
applying the ALS iterative scheme in the canonical format, see e.g. \cite{KoldaB:08}.)

{\bf 2. The  Hartree-Fock equation in AO basis set.} 

The $2N$-electrons Hartree-Fock equation for pairwise $L^2$-orthogonal
electronic orbitals, $\psi_i: \mathbb{R}^3\to\mathbb{R} $,
$\psi_i\in H^1(\mathbb{R}^3)$, reads as
\begin{equation} \label{HaFo eq.}
  {\cal F}  \psi_i({ x}) = \lambda_i \, \psi_i({ x}),\quad
  \int_{\mathbb{R}^3}\psi_i\psi_j  dx =\delta_{ij},
   \; i,j=1,...,N_{orb}
\end{equation}
where the nonlinear  Fock operator ${\cal F}$ is given by
\[
{\cal F} :=-\frac{1}{2} \Delta + V_c(\cdot) + V_H(\cdot) + {\cal K}.
\]
Here the nuclear potential takes the form  $V_c(x)=- \sum_{\nu=1}^{M}\frac{Z_\nu}{\|{x} -a_\nu \|}$, $Z_\nu >0$, $a_\nu\in \mathbb{R}^3$,
while the Hartree potential $V_H({ x})$ and the nonlocal exchange operator $\cal K$ read as
 \begin{equation}
   V_{H}({x}):= \rho \star \frac{1}{\|\cdot\|} 
=  \int_{\mathbb{R}^3}
 \frac{\rho({y})}{\|{x}-{y}\|}\, d{y} , \quad{x}\in {\mathbb{R}^3},
\end{equation} 
and 
\begin{equation}
\left( {\cal K}\psi\right) ({ x})
:=-  \frac{1}{2}
\sum_{i=1}^{N_{orb}}\left(\psi \, \psi_i\star \frac{1}{\|\cdot\|} \right)
\psi_i (x)
=  - \frac{1}{2}\int_{\mathbb{R}^3} \; 
\frac{\tau({ x}, { y})}{\|{x} - { y}\|}\,\psi({y}) d{ y},
\label{VExch}
\end{equation} 
respectively. Conventionally, we use the definitions
\begin{equation*} \label{DMatr-eq}
\tau({ x}, { y}) := 2\sum_{i=1}^{N_{orb}} \; \psi_i({ x}) \psi_i({ y}), 
\quad  \rho( {x}) :=\tau({x},{x}),\quad
\end{equation*}
for the density matrix $\tau({x}, { y})$, and electron density $\rho({x})$.

Usually, the Hartree-Fock equation is approximated by the standard Galerkin projection of the 
initial problem  (\ref{HaFo eq.}) by using the physically justified reduced basis sets 
(say, GTO type orbitals).
For a given finite Galerkin basis set $\{g_\mu \}_{1\leq \mu \leq N_b}$,
$g_\mu\in H^1(\mathbb{R}^3) $, the occupied
molecular orbitals $\psi_i$ are represented (approximately) as
$
\psi_i=\sum\limits_{\mu=1}^{N_b} C_{\mu i} g_\mu, \quad i=1,...,N_{orb}.
$
To derive an equation for the unknown coefficients matrix
$C=\{C_{ \mu i} \}\in \mathbb{R}^{N_b \times  N_{orb}}$, first, we introduce
the mass (overlap) matrix $S=\{S_{\mu \nu} \}_{1\leq \mu, \nu \leq N_b}$, given by
$
S_{\mu \nu}=\int_{\mathbb{R}^3} g_\mu g_\nu  dx,
$
and the stiffness matrix $H=\{h_{\mu \nu}\}$ of the core Hamiltonian
${\cal H}=-\frac{1}{2} \Delta + V_c $ (the single-electron integrals),
\[
h_{\mu \nu}= \frac{1}{2} \int_{\mathbb{R}^3}\nabla g_\mu \cdot \nabla g_\nu dx +
\int_{\mathbb{R}^3} V_c(x) g_\mu g_\nu dx, \quad 1\leq \mu, \nu \leq N_b.
\]
The core Hamiltonian matrix $H$ can be precomputed in $O(N_b^2)$ operations via grid-based approach.

Given the finite basis set $\{g_\mu \}_{1\leq \mu \leq N_b}$, $g_\mu\in H^1(\mathbb{R}^3) $, the
associated fourth order two-electron integrals (TEI) tensor,
${\textbf{B}}=[b_{\mu \nu \lambda \sigma}]$, is defined entrywise by (\ref{eqn:TEI_tens}), 
where $\mu, \nu, \lambda, \sigma \in \{1,...,N_b\}=:{\cal I}_b$.
In computational quantum chemistry the nonlinear terms representing the Galerkin 
approximation to the Hartree and
exchange operators are calculated traditionally by using the low-rank Cholesky decomposition 
of a matrix associated with the
TEI tensor ${\bf B}=[b_{\mu \nu \kappa \lambda}]$ defined in (\ref{eqn:TEI_tens}), 
that initially has the computational and storage complexity  of order $O(N_b^4)$.
 
Introducing the $N_b\times N_b$ matrices $J(D)$ and $K(D)$,
\begin{equation*}\label{C-K matr}
J(D)_{\mu \nu}= \sum\limits_{\kappa, \lambda=1}^{N_b}
b_{\mu \nu, \kappa \lambda}D_{\kappa \lambda},\quad
K(D)_{\mu \nu}= - \frac{1}{2} \sum\limits_{\kappa, \lambda=1}^{N_b}
b_{\mu \lambda, \nu \kappa}D_{\kappa \lambda},
\end{equation*}
where  $D= 2 C C^T \in \mathbb{R}^{N_b\times N_b}$ is the rank-$N_{orb}$ symmetric density matrix,
one then represents the complete Fock matrix $F$ by
\begin{equation*}\label{Fock_matr}
 F(D)= H+ J(D) + K(D).   
\end{equation*}
The resultant Galerkin system of nonlinear equations
for the coefficients matrix $C\in \mathbb{R}^{N_b \times N_{orb}}$, 
and the respective eigenvalues $\Lambda= diag(\lambda_1,...,\lambda_{N_b})$, reads as
\begin{align} \label{HF discr}
  F(D) C &= SC \Lambda, \quad C^T SC =  I_N,   \nonumber
\end{align}
where the second equation represents the orthogonality constraints
$\int_{\mathbb{R}^3}\psi_i\psi_j dx=\delta_{ij}$, and
$I_N$ denotes the $N_b\times N_b $ identity matrix.\\

{\bf 3. Tensor approximation to the Newton kernel in 3D.} 

Methods of separable approximation to multivariate spherically symmetric functions
by using the Gaussian sums have been addressed in the chemical and mathematical literature 
since \cite{Boys:56} and \cite{Braess:95,Stenger}, respectively.

In this section, we discuss for the readers convenience the grid-based method for 
the low-rank canonical and Tucker tensor 
representations of a spherically symmetric functions $p(\|x\|)$, 
$x\in \mathbb{R}^d$ in the particular case of the 3D Newton kernel
$p(\|x\|)=\frac{1}{\|x\|}$, $x\in \mathbb{R}^3$
by its projection onto the set
of piecewise constant basis functions, see \cite{BeHaKh:08} for more details.

In the computational domain  $\Omega=[-b/2,b/2]^3$, 
let us introduce the uniform $n \times n \times n$ rectangular Cartesian grid $\Omega_{n}$
with the mesh size $h=b/n$ (usually, $n=2k$).
Let $\{ \psi_\textbf{i}\}$ be a set of tensor-product piecewise constant basis functions,
$  \psi_\textbf{i}(\textbf{x})=\prod_{\ell=1}^d \psi_{i_\ell}^{(\ell)}(x_\ell)$,
for the $3$-tuple index ${\bf i}=(i_1,i_2,i_3)$, $i_\ell \in \{1,...,n\}$, $\ell=1,\, 2,\, 3 $.
The kernel $p(\|x\|)$ can be discretized by its projection onto the basis set $\{ \psi_\textbf{i}\}$
in the form of a third order tensor of size $n\times n \times n$, defined entrywise as
\begin{eqnarray}
\mathbf{P}:=[p_{\bf i}] \in \mathbb{R}^{n\times n \times n},  \quad
 p_{\bf i} = 
\int_{\mathbb{R}^3} {\psi_{{\bf i}}({x})}p({\|{x}\|}) \,\, \mathrm{d}{x}.
  \label{galten}
\end{eqnarray}

Given $M$, the low-rank canonical decomposition of the $3$rd order tensor $\mathbf{P}$ is based 
on using exponentially convergent 
$\operatorname*{sinc}$-quadratures for approximation of the Laplace-Gauss transform 
to the analytic function $p(z)= 1/z$ as follows,
\begin{equation} \label{eqn:sinc_Newt}
 1/{\|{x}\|} =  \frac{2}{\sqrt{\pi}} \int_{\mathbb{R}_+} 
e^{- t^2\|{x}\|^2} \,\mathrm{d}t 
\approx 
\sum_{k=-M}^{M} a_k e^{- t_k^2\|{x}\|^2}= 
\sum_{k=-M}^{M} a_k  \prod_{\ell=1}^3 e^{-t_k^2 x_\ell^2},
\end{equation}
where the quadrature points and weights are given by 
\begin{equation*} \label{eqn:hM}
t_k=k \mathfrak{h}_M , \quad a_k= \frac{2}{\sqrt{\pi}} \mathfrak{h}_M, \quad 
\mathfrak{h}_M=C_0 \log(M)/M , \quad C_0>0.
\end{equation*}

Under the assumption $0< a \leq \|x \|  < \infty$, $x=(x_1,x_2,x_3)\in \mathbb{R}^3$,
this quadrature can be proven to provide the exponential convergence rate in $M$
for a class of analytic functions, see \cite{Stenger,Braess:95,HaKhtens:04I},
 \begin{equation*} \label{sinc_conv}
\left|1/{\|{x}\|} - \sum_{k=-M}^{M} a_k e^{- t_k^2\|{x}\|^2} \right|  
\le \frac{C}{a}\, \displaystyle{e}^{-\beta \sqrt{M}},  
\quad \text{with some} \ C,\beta >0.
\end{equation*}
Combining \eqref{galten} and \eqref{eqn:sinc_Newt}, and taking into account the 
separability of the Gaussian basis functions, we arrive at the low-rank 
approximation to each entry of the tensor $\mathbf{P}$,
\begin{equation*} \label{eqn:C_nD_0}
 p_{\bf i} \approx \sum_{k=-M}^{M} a_k   \int_{\mathbb{R}^3}
 \psi_{\bf i}(\textbf{x}) e^{- t_k^2\|\mathbf{x}\|^2} \mathrm{d}{\bf x}
=  \sum_{k=-M}^{M} a_k  \prod_{\ell=1}^{3}  \int_{\mathbb{R}}
\psi^{(\ell)}_{i_\ell}(x_\ell) e^{- t_k^2 x^2_\ell } \mathrm{d}x_\ell.
\end{equation*}
Define the vector (recall that $a_k >0$)
\begin{equation*} \label{eqn:galten_int}
\textbf{p}^{(\ell)}_k
= a_k^{1/3} \left[b^{(\ell)}_{i_\ell}(t_k)\right]_{i_\ell=1}^{n} \in \mathbb{R}^{n}
\mbox { with }
  b^{(\ell)}_{i_\ell}(t_k)= 
\int_{\mathbb{R}} \psi^{(\ell)}_{i_\ell}(x_\ell) e^{- t_k^2 x^2_\ell } \mathrm{d}x_\ell,
\end{equation*}
then the $3$rd order tensor $\mathbf{P}$ can be approximated by 
the $R$-term  canonical representation
\begin{equation} \label{eqn:sinc_general}
    \mathbf{P} \approx  \mathbf{P}_R =
\sum_{k=-M}^{M} a_k \bigotimes_{\ell=1}^{3}  {\bf b}^{(\ell)}(t_k)
= \sum\limits_{q=1}^{R} {\bf p}^{(1)}_q \otimes {\bf p}^{(2)}_q \otimes {\bf p}^{(3)}_q
\in \mathbb{R}^{n\times n \times n},
\end{equation}
where $R=2M+1$, and the canonical vectors are renumbered by $k \to q=k+M+1$, 
${\bf p}^{(\ell)}_q ={\bf p}^{(\ell)}_{k} \in \mathbb{R}^n$, $\ell=1, 2, 3$.
For the given threshold $\varepsilon >0 $,  $M=O(|\log \varepsilon|^2)$ is chosen as the minimal number,
such that in the max-norm we have
\begin{equation*} \label{eqn:error_control}
\| \mathbf{P} - \mathbf{P}_R \|  \le \varepsilon \| \mathbf{P}\|.
\end{equation*}
The symmetric canonical tensor ${\bf P}_{R}$ in (\ref{eqn:sinc_general})
approximates the discretized 3D symmetric kernel function 
$p({\|x\|})= 1/\|x\|$ ($x\in \Omega$), 
centered at the origin, such that ${\bf p}^{(1)}_q={\bf p}^{(2)}_q={\bf p}^{(3)}_q$ ($q=1,...,R$).

\footnotesize{

 }

\end{document}